\documentclass[10pt]{amsart}

\newcommand\dss{\displaystyle}

\usepackage{latexsym, amsfonts, amssymb}

\newenvironment{thm}{\subsection{}{\textbf {Theorem.}}\em}{}
\newenvironment{prop}{\subsection{}{\textbf {Proposition.}}\em}{}
\newenvironment{cor}{\subsection{}{\textbf {Corollary.}}\em}{}
\newenvironment{lem}{\subsection{}{\textbf {Lemma.}}\em}{}

\newenvironment{pf}{\noindent{\textbf {Proof.}}}
{\begin{flushright}\eop \end{flushright}\smallskip}

\newenvironment{defn}{\subsection{}{\textbf {Definition.}}\em}{\smallskip}
\newenvironment{eg}{\subsection{}{\textbf {Example.}}}{\smallskip}

\newenvironment{rem}{\subsection{}{\textbf {Remark.}}}{\smallskip}

\newcommand\cB{\ensuremath{\mathcal B}}
\newcommand\cC{\ensuremath{\mathcal C}}

\newcommand\cE{\ensuremath{\mathcal E}}

\newcommand\cH{\ensuremath{\mathcal H}}

\newcommand\cJ{\ensuremath{\mathcal J}}
\newcommand\cK{\ensuremath{\mathcal K}}

\newcommand\cM{\ensuremath{\mathcal M}}
\newcommand\cN{\ensuremath{\mathcal N}}

\newcommand\cU{\ensuremath{\mathcal U}}

\newcommand\bbC{\ensuremath{\mathbb C}}

\newcommand\bbM{\ensuremath{\mathbb M}}
\newcommand\bbN{\ensuremath{\mathbb N}}

\newcommand\bbR{\ensuremath{\mathbb R}}

\newcommand\bbT{\ensuremath{\mathbb T}}

\newcommand\bbZ{\ensuremath{\mathbb Z}}

\newcommand\hilb{\ensuremath{\mathcal H}}

\newcommand\ol{\ensuremath{\overline}}

\newcommand\eop{{{\hfil \ensuremath \Box}}}
\newcommand\eps{\ensuremath {\varepsilon}}
\newcommand\norm{\ensuremath {\Vert}}

\newcommand\bofh{\ensuremath{\cB ( \cH)}}

\newcommand\kofh{\ensuremath{\cK ( \cH)}}

\newcommand\nul{\ensuremath {\mathrm {nul}}}

\usepackage{mathrsfs}

\usepackage{amssymb, amsmath}
\usepackage{bm}

\usepackage[normalem]{ulem}

\renewcommand{\phi}{\varphi}
\renewcommand{\rho}{\varrho}

\newcommand{\Norm}[1]{\ensuremath{\left\lVert {#1} \right\rVert}}

\newcommand{\Field}[1]{\ensuremath{\mathbb{#1}}}
\newcommand{\N}{\Field{N}}

\newcommand{\C}{\Field{C}}

\newcommand{\R}{\Field{R}}

\newcommand{\AND}{\ \ \ \text{and}\ \ \ }

\newcommand{\Mn}[1]{\ensuremath{\Field{M}_{#1}}}

\usepackage{xcolor}
\definecolor{vvvlgrey}{rgb}{.92,.92,.92}
\definecolor{vvlgrey}{rgb}{.9,.9,.9}
\definecolor{vlgrey}{rgb}{.8,.8,.8}
\definecolor{lgrey}{rgb}{.5,.5,.5}
\definecolor{darkgreen}{rgb}{0.3,.47,0.33}
\definecolor{purple}{rgb}{.78,0,0.68}
\definecolor{dullmagenta}{rgb}{0.4,0,0.4}   
\definecolor{darkblue}{rgb}{0,0,0.4}

\newcommand{\offdiag}{property \textsc{(CN)}}
\newcommand{\offdr}{property \textsc{(CR)}}

\begin{document}
\title{Hilbert space operators with compatible off-diagonal corners}

\thanks{${}^1$ Research supported in part by NSERC (Canada)}
\thanks{{\ifcase\month\or Jan.\or Feb.\or March\or April\or May\or 
June\or
July\or Aug.\or Sept.\or Oct.\or Nov.\or Dec.\fi\space \number\day,
\number\year}}
\author
	[L. Livshits]{{L. Livshits}}
\address
	{Department of  Mathematics\\
	Colby College\\
	Waterville, ME\\
	USA\ \ \ 188745}
\email{llivshi@colby.edu}
\author
	[G. MacDonald]{{G. MacDonald${}^1$}}
\address
	{School of Mathematical and Computational Sciences\\
	University of Prince Edward Island\\
	Charlottetown, PE \\
	Canada  \ \ \ C1A~4P3}
\email{gmacdonald@upei.ca}
\author
	[L.W. Marcoux]{{L.W.~Marcoux${}^1$}}
\address
	{Department of Pure Mathematics\\
	University of Waterloo\\
	Waterloo, ON \\
	Canada  \ \ \ N2L~3G1}
\email{LWMarcoux@uwaterloo.ca}
\author
	[H. Radjavi]{{H. Radjavi${}^1$}}
\address
	{Department of Pure Mathematics\\
	University of Waterloo\\
	Waterloo, ON \\
	Canada  \ \ \ N2L~3G1}
\email{hradjavi@uwaterloo.ca}

\begin{abstract}
Given a complex, separable Hilbert space $\hilb$, we characterize those operators for which $\norm P T (I-P) \norm = \norm (I-P) T P \norm$ for all orthogonal projections $P$ on $\hilb$.   When $\hilb$ is finite-dimensional, we also obtain a complete characterization of those operators for which $\mathrm{rank}\, (I-P) T P = \mathrm{rank}\, P T (I-P)$ for all orthogonal projections $P$.  When $\hilb$ is infinite-dimensional, we show that any operator with the latter property is normal, and its spectrum is contained in either a line or a circle in the complex plane.   
\end{abstract}

\maketitle
\markboth{\textsc{  }}{\textsc{}}


 
\section{Introduction} \label{sec1}


\subsection{} \label{sec1.01}

Let $\hilb$ be a complex, separable Hilbert space.     By $\bofh$, we denote the algebra of bounded, linear operators on $\hilb$.   If $\dim\, \hilb = n < \infty$, then we identify $\hilb$ with $\bbC^n$ and $\bofh$ with $\bbM_n(\bbC)$.

One of the most important open problems in operator theory is the \emph{Invariant Subspace Problem}, which asks whether or not every bounded, linear operator $T$ acting on a complex, infinite-dimensional, separable Hilbert space $\hilb$ admits a non-trivial invariant subspace; that is, a closed subspace $\cM \not \in \{ \{ 0 \}, \hilb\}$ for which $T \cM \subseteq \cM$. 

We say that an operator $T \in \bofh$ is \textbf{(orthogonally) reductive} if for each orthogonal projection $P \in \bofh$, the condition $P T (I-P) = 0$ implies that $(I-P) T P = 0$.  The \emph{Reductive Operator Conjecture} is the assertion that every reductive operator is  normal.  It was shown by Dyer, Pederson and Procelli~\cite{DPP1972} that the Invariant Subspace Problem admits a positive solution if and only if the \emph{Reductive Operator Conjecture} is true.    

Our goal in this paper is to study two variants of orthogonal reductivity.  Let $T \in \bofh$ and $P \in \bofh$ be an orthogonal projection.  We refer to the operator 
\[
P^\perp T P: P\hilb \to P^\perp \hilb\]
as an \textbf{off-diagonal corner} of $T$.

Relative to the decomposition $\hilb = P \hilb \oplus P^\perp \hilb$, we may write $T=\left[\begin{smallmatrix}
 A & B \\ C & D 
\end{smallmatrix}\right]  $.  We refer to the block-entries of such block-matrices via their geographic positions: NW, NE, SE, SW, and the NE and the SW block-entries are examples of the {{off-diagonal corners}}.

In the work below, we shall be interested in two phenomena:  firstly, when the operator norm of $B (= B_P)$ coincides with the operator norm of $C (=C_P)$ for all projections $P$, and secondly, when the rank of $B$ coincides with the rank of $C$ for all projections $P$.  Clearly, any operator which satisfies one of these two conditions is orthogonally reductive.  An example is given in Section~\ref{sec5} below to show that the converse to this statement is false.
 
In the case of normal matrices, some related work has been done by Bhatia and Choi~\cite{BC2006}.   For instance, if the dimension of the space is $2 n < \infty$, and if $P$ is a projection of rank $n$, it is a consequence of the fact that the Euclidean norm of the $k^{th}$ column of a normal matrix coincides with that of the $k^{th}$ row for all $k$ that the Hilbert-Schmidt (or Frobenius) norm of $B$ always equals that of $C$.  Further, they show that $\norm B \norm \le \sqrt{n} \norm C \norm$, and that equality can be achieved for some normal matrix $T \in \bbM_{2n}(\bbC)$ and some projection $P$ of rank $n$ if and only if $n \le 3$. 


\begin{defn} \label{defn1.02}
Let $T \in \bofh$.   We say that $T$ has the \textbf{common norm property} \emph{(}{\offdiag}\emph{)}  if for any projection $P \in \bofh$ we have that 
\[
\norm P T P^\perp \norm = \norm P^\perp T P \norm. \]

We denote by $\mathfrak{G_{norm}}$ the set of operators with {\offdiag}.
We say that $T$ has the \textbf{common rank property} \emph{(}{\offdr}\emph{)} if for any projection $P \in \bofh$ we have that 
\[
\mathrm{rank}\, P T P^\perp = \mathrm{rank}\, P^\perp T P. \]

We denote by $\mathfrak{G_{rank}}$ the set of operators with {\offdr}.
\end{defn}

As we shall see, our results depend upon whether or not $\hilb$ is finite-dimensional.  When the Hilbert space is finite-dimensional and of dimension at least four, then we shall show that the set of operators satisfying {\offdiag} coincides with the set of operators satisfying {\offdr}, and that this consists of those operators which are scalar translates of scalar multiples of hermitian  (or of unitary) operators.  (See Theorem~\ref{thm3.16} below.)

In the infinite-dimensional setting, we obtain a complete characterization of those operators satisfying {\offdiag}.  Again, any scalar translate of a scalar multiple of a hermitian operator will suffice.  This time, however, the unitary operators involved must have essential spectrum contained in only half of a circle.  (See Theorem~\ref{thm4.13} below.)

The problem of characterizing those  operators acting on an infinite-dimensional Hilbert space which enjoy {\offdr} is much more delicate.    We are able to demonstrate that any operator $T$ satisfying {\offdr} must once again be a scalar translate of a scalar multiple of a hermitian (or of a unitary) operator.  In particular, such operators are normal.   However, an obstruction occurs in that it is not the case that every unitary operator has  {\offdr}.  
Indeed, as is well-known (see Section~\ref{sec5} for an example) -- not every unitary operator is reductive. 



\subsection{} \label{sec1.03}
We shall need some standard notation and definitions in what follows.   

If $T=\left[\begin{smallmatrix}
 A & B \\ C & D 
\end{smallmatrix}\right]$ is a block-matrix in $\bbM_n(\bbC)$, and $A$ is invertible, then the matrix $D-CA^{-1}B$ is said to be the {\textbf{Schur complement}} of $A$ in $T$ and is denoted by $T|A$. In such a case $T$ is invertible if and only if $T|A$ is, and when this happens, the SE block-corner $\left(T^{-1}\right)_{_{SE}}$ of $T^{-1}$ is $(T|A)^{-1}$. Furthermore:
\[%
\left(T^{-1}\right)_{_{SW}}=-(T|A)^{-1}CA^{-1}\ \text{ and }\ \left(T^{-1}\right)_{_{NE}}=-A^{-1}B\,(T|A)^{-1}.
\]%

Similarly, if $B$ is invertible then $C-DB^{-1}A$ is the Schur complement $T|B$ of $B$ in $T$, and $T$ is invertible if and only if $T|B$ is, in which case $$\left(T^{-1}\right)_{_{NE}}=(T|B)^{-1}.$$
Corresponding statements and concepts apply to $C$ and $D$ as well.

As always, $\bbT = \{ z \in \bbC: |z| = 1\}$. A subset of $\C$ is \textbf{circlinear} if it is contained in a circle or a straight line. By $\kofh$, we denote the closed, two-sided ideal of compact operators in $\bofh$, and $\pi: \bofh \to \bofh/\kofh$ denotes the canonical map from $\bofh$ into the \textbf{Calkin algebra} $\bofh/\kofh$.
The \textbf{essential spectrum} $\sigma_e(T)$ of $T \in \bofh$ is the spectrum of $\pi(T)$ in the Calkin algebra $\bofh/\kofh$, and we say that $T$ is a \textbf{Fredholm operator} if $0 \not \in \sigma_e(T)$.   The \textbf{Fredholm domain} of $T$ is $\rho_{F}(T) = \bbC \setminus \sigma_e(T)$.   We say that  $T$ is a \textbf{semi-Fredholm operator} if $\pi(T)$ is either left or right invertible in $\bofh/\kofh$, and define the \textbf{semi-Fredholm domain}  of $T$ to be $\rho_{sF}(T) = \{ \lambda \in \bbC:  (T-\lambda I) \mbox{ is semi-Fredholm}\}$.  The complement of $\rho_{sF}(T)$ is called the \textbf{left-right essential spectrum} of $T$ and is denoted by $\sigma_{\ell r e}(T)$.  If $T$ is semi-Fredholm, we define the \textbf{index} of $T$ to be $\mathrm{ind}\, T = \nul\, T - \nul\, T^* \in \bbZ \cup \{-\infty, \infty\}$.  When $T$ is Fredholm, we have that $\mathrm{ind}\, T \in \bbZ$.

We say that $T$ is \textbf{triangular} if there exists an orthonormal basis $\{ e_n\}_{n=1}^\infty$ for $\hilb$ such that the matrix $[T] = [t_{i, j}]$ for $T$ relative to this basis (i.e. $t_{i, j} = \langle T e_j, e_i \rangle$) satisfies $t_{i, j} = 0$ for all $i > j$.  The operator is said to be \textbf{quasitriangular} if it is of the form $T = T_0 + K$, where $T_0$ is triangular and $K$ is compact.   It was shown by Apostol, Foia\c{s}, and Voiculescu~\cite{AFV1973} that $T$ is quasitriangular if and only if $\mathrm{ind}\, (T-\lambda I) \ge 0$ whenever $\lambda \in \rho_{sF}(T)$.  Finally, $T$ is \textbf{biquasitriangular} if each of $T$ and $T^*$ is quasitriangular, i.e. if and only if $\mathrm{ind}\, (T-\lambda I) = 0$ for all $\lambda \in \rho_{sF}(T)$.

Recall also that if $T \in \bofh$, then $|T| = (T^* T)^{1/2}$ denotes the absolute value of $T$.  A unitary operator $U \in \bofh$ is said to be \textbf{absolutely continuous} if the spectral measure for $U$ is absolutely continuous with respect to Lebesgue measure restricted to $\sigma(U)$, while $U$ is said to be \textbf{singular} if the spectral measure of $U$ is singular with respect to Lebesgue measure restricted to $\sigma(U)$.  These notions will only be used in Section~\ref{sec5}.

\section{Preliminary results} \label{sec2}


\subsection{} \label{sec2.01}
We begin with a few simple remarks.  Although the proofs are rather elementary, we shall list these in the form of a Proposition so as to be able to more easily refer to them later.  The proofs are left to the reader.


\begin{prop} \label{prop2.02}
Suppose that $R, T \in \bofh$ and that $R$ has {\offdr} and $T$ has {\offdiag}.
\begin{enumerate}  
	\item[(a)] For all $\lambda, \mu \in \bbC$, we have that
	\begin{itemize}
		\item{}
		$\lambda I + \mu R$ and $R^*$ have {\offdr}, while 
		\item{}
		$\lambda I + \mu T$ and $T^*$ have {\offdiag}.   
	\end{itemize}
	\item[(b)]
	Suppose that $\hilb = \hilb_1 \oplus \hilb_2$. 
	\begin{itemize}
		\item{}
		If there exist $A \in \cB(\hilb_1)$ and $D \in \cB(\hilb_2)$ such that $R = A \oplus D$, then $A$ and $D$ both have {\offdr}.
		\item{}
		If there exist $A \in \cB(\hilb_1)$ and $D \in \cB(\hilb_2)$ such that $T = A \oplus D$, then $A$ and $D$ both have {\offdiag}.
	\end{itemize}	
	\item[(c)]
	If $V \in \bofh$ is unitary, then 
	\begin{itemize}
		\item{}
		$V^* R V$ has {\offdr} and 
		\item{}
		$V^* T V$ has {\offdiag}.
	\end{itemize}	
	\item[(d)]
	If $L = L^* \in \bofh$, then $L$ has both {\offdr} and {\offdiag}.
\end{enumerate}
\end{prop}

%


\bigskip

In the case of {\offdiag}, we also observe the following.  For $T \in \bofh$, let us denote by $\cU(T)$ the \textbf{unitary orbit} of $T$, i.e. $\cU(T) = \{ V^* T V: V \in \bofh \mbox{ unitary} \}$.  Recall that two operators $S$ and $T$ are said to be \textbf{approximately unitarily equivalent} if $S \in \ol{\cU(T)}$ (equivalently, $T \in \ol{\cU(S)}$). The proofs of the following assertions are elementary and are left to the reader.

\smallskip

\begin{prop} \label{prop2.03}  
\begin{enumerate}
	\item[(a)]
	The set $\mathfrak{G_{norm}}$ of operators with {\offdiag} is closed.
	\item[(b)]
	If $T \in \bofh$ has {\offdiag} and there exists $S \in \ol{\cU(T)}$ of the form $S = A \oplus D$, then $A, D$ have {\offdiag}.
\end{enumerate}
\end{prop}

%


\bigskip

The following remark, while innocuous in appearance, is actually the key to a number of calculations below.

\smallskip

\begin{rem} \label{rem2.04}
Let $U \in \bofh$ be a unitary operator and $P \in \bofh$ be a projection.  Write $U = \left[\begin{smallmatrix}
 A & B \\ C & D 
\end{smallmatrix}\right]$ relative to $\hilb = P \hilb \oplus P^\perp \hilb$.    
	The fact that $U$ is unitary implies that 
	\[
	I = A A^* + B B^* = A^* A + C^*C. \]
	Thus $B B^* = I - A A^*$ and $C^* C = I - A^* A.$  
	
	It follows that 
	\[
	\norm B \norm^2 = \norm B B^* \norm = 1 - \min \{ \lambda: \lambda \in \sigma(A A^*) \}, \]
	and similarly 
	\[
	\norm C \norm^2 = \norm C^* C \norm = 1 - \min \{ \mu : \mu \in \sigma(A^*A)\}. \]
	 
	However, it is a standard fact that $\sigma(A A^*) \cup \{ 0 \} = \sigma (A^* A) \cup \{ 0\}$, and thus the only way that we can have $\norm B \norm \ne \norm C \norm$ is if either 
	\begin{enumerate}
		\item[(I)]
		$0 \in \sigma(A A^*)$ but $0 \not \in \sigma(A^*A)$, or 
		\item[(II)]
		$0 \in \sigma(A^* A)$ but $0 \not \in \sigma(A A^*)$.
	\end{enumerate}
		
	This argument demonstrates the rather interesting fact that if $U = \left[\begin{smallmatrix}
 A & B \\ C & D 
\end{smallmatrix}\right]$ is a unitary operator and $\norm B \norm \ne \norm C \norm$, then 
	\[
	\min(\norm B \norm, \norm C \norm) < 1 = \max(\norm B \norm, \norm C \norm).\]
	
	{In particular, if $U = \left[\begin{smallmatrix}
 A & B \\ C & D 
\end{smallmatrix}\right]$ is a unitary such that $\norm B \norm < \norm C \norm (=1)$, then every unitary $U'$ close enough to $U$ has the form $\left[\begin{smallmatrix}
 A' & B' \\ C' & D' 
\end{smallmatrix}\right]$ where $\norm B' \norm < \norm C' \norm \bm{=1},$ which is remarkable.}

\end{rem}

\section{The finite-dimensional setting} \label{sec3}

\subsection{} \label{sec3.01}
We now turn to the case where the Hilbert space under consideration is finite-dimensional (and complex).   


\begin{prop} \label{prop3.02}
Let $n \ge 2$ be an integer and $T \in \bbM_n(\bbC)$.   If $T$ has {\offdiag} or {\offdr}, then $T$ is normal.
\end{prop}

\begin{pf}
This is an easy consequence of the fact that given any $T \in \bbM_n(\bbC)$, there exists an orthonormal basis with respect to which the matrix of $T$ is upper triangular.   Either property clearly implies that the matrix of $T$ is in fact diagonal with respect to this basis, and hence that $T$ is normal.
\end{pf}


\begin{prop} \label{prop3.03}
Let $n \ge 2$, and let $U \in \bbM_n(\bbC)$ be a unitary operator.   Then $U$ has both {\offdiag} and {\offdr}.
\end{prop}

\begin{pf}
Let $P \in \bbM_n(\bbC)$ be a projection, and relative to the decomposition $\bbC^n = P \bbC^n \oplus P^\perp \bbC^n$, let us write 
\[
U = \begin{bmatrix} A &B \\ C & D \end{bmatrix}. \]
As noted in Remark~\ref{rem2.04}, since $U$ is unitary, we have that $B B^* = I - A A^*$ and $C^* C = I - A^*A$.  Observe, however, that in the finite-dimensional setting we have that $A^* A$ is unitarily equivalent to $A A^*$, and thus $B B^*$ is unitarily equivalent to $C^* C$.   Thus 
\begin{itemize}
	\item{}
	$\norm B \norm = \norm C \norm$, and 
	\item{}
	$\mathrm{rank}\, B = \mathrm{rank}\, B B^* = \mathrm{rank}\, C^* C = \mathrm{rank}\, C$.
\end{itemize}	
\end{pf}


Combining this with Proposition~\ref{prop2.02} (a), we obtain:

\begin{prop} \label{prop3.04}
Let $T \in \bbM_n(\bbC)$.  If $T$ is either hermitian or unitary, then for all $\lambda, \mu \in \bbC$, $\lambda I + \mu T$ has both {\offdiag} and {\offdr}.
\end{prop}

\smallskip

Our goal is to prove that if $T \in \bbM_n(\bbC)$ has either {\offdiag} or {\offdr}, then it is of the form $\lambda I + \mu X$ where $X$ is either hermitian or unitary.


\begin{rem} \label{rem3.05}
The common link between these two cases is the geometry of the set of eigenvalues of $T$.   If $T$ is normal, then $T = \lambda I + \mu V$ where $V$ is unitary if and only if all of the eigenvalues of $T$ lie on a common circle.   If $T$ is normal, then $T = \lambda I + \mu L$ where $L = L^*$ if and only if all of the eigenvalues of $T$ lie on a common line.  That is to say, the union of these two sets of operators is precisely the class of normal operators whose spectra are circlinear.    

Given a matrix $B \in \bbM_{n, m}(\bbC)$, we denote by $\norm B \norm_2 = \mathrm{tr} (B^* B)^{1/2}$ the Hilbert-Schmidt (or Fr\"{o}benius) norm of $B$.
\end{rem}


\begin{prop} \label{prop3.06}
Let $k, \ell \ge 1$ be integers, and suppose that $T = \begin{bmatrix} A & B \\ C & D \end{bmatrix}$ is a normal operator in $\cB(\bbC^k \oplus \bbC^\ell)$.   Then 
\[
\norm B \norm_2 = \norm C \norm_2. \]
\end{prop}

\begin{pf}
The fact that $T$ is normal implies that $A A^* + B B^*  = A^* A + C^* C$.  Using the fact that the trace is linear and that $\mathrm{tr} (X Y) = \mathrm{tr} (Y X)$ for all $X \in \bbM_{k, \ell}(\bbC)$,  $Y \in \bbM_{\ell, k}(\bbC)$, we see that 
\[
\norm B \norm_2 = \mathrm{tr}(B B^*) = \mathrm{tr}(C^* C) = \norm C \norm_2. \]
\end{pf}


We begin by considering the exceptional cases where the dimension of the underlying Hilbert space is too small to allow anything interesting to happen.

\smallskip

\begin{prop} \label{prop3.07}
If $2 \le n \le 3$, and let $T \in \bbM_n(\bbC)$.  The following are equivalent:
\begin{enumerate}
	\item[(a)]
	$T$ is normal.
	\item[(b)]
	$T$ has {\offdiag}.
	\item[(c)]
	$T$ has {\offdr}.
\end{enumerate}	
\end{prop}

\begin{pf}
By Proposition~\ref{prop3.02}, both (b) and (c) imply (a).

Conversely, if $T \in \bbM_n(\bbC)$ is normal and $0 \ne P \ne I$ is a projection in $\bbM_n(\bbC)$, then $P T P^\perp$ and $P^\perp T P$ both have rank at most one.     From this and from Proposition~\ref{prop3.06}, we find that 
\[
\norm P T P^\perp \norm = \norm P T P^\perp\norm_2 = \norm P^\perp T P \norm_2  = \norm P^\perp T P \norm. \]

Thus $T$ has {\offdiag}; that is, (a) implies (b).  This also shows that $P T P^\perp$ and $P^\perp T P$ either both have rank $0$ or both have rank $1$.  Hence (a) implies (c) as well.
\end{pf}

For the remainder of this section, we shall assume that the dimension $n$ of the underlying Hilbert space is at least $4$.


\begin{rem} \label{rem3.08}
Let us now show that the problem of characterizing which operators in $\bbM_n(\bbC)$ have {\offdiag} (resp. {\offdr}) reduces to the case where $n= 4$.  Of course, by Proposition~\ref{prop3.02}, we may restrict our attention to normal operators.

Let $n > 4$, and suppose that $T \in \bbM_n(\bbC)$ is normal.

\begin{itemize}
	\item
	As observed in Remark~\ref{rem3.05}, if all of the eigenvalues of $T$ are either co-linear or co-circular (i.e. all lie on the same circle), then there exist $\alpha, \beta \in \bbC$ and either a hermitian operator $L$ or a unitary operator $V$ such that $T = \alpha I + \beta L$, or $T = \alpha I + \beta V$.  Either way, by Proposition~\ref{prop3.04}, $T$ has {\offdiag} and {\offdr}.\smallskip
	\item
	Conversely, suppose that $T$ has {\offdiag} (resp. $T$ has {\offdr}), and \emph{suppose} we know that every $X \in \bbM_4(\bbC)$ with {\offdiag} (resp. with {\offdr}) has eigenvalues that are either co-linear or co-circular.   Given any $\{ \theta_1, \theta_2, \theta_3, \theta_4\}$ in $\sigma(T)$, we can write $T = R \oplus Y$, where $R$ is a normal operator in $\bbM_4(\bbC)$ with $\sigma(R) = \{\theta_1, \theta_2, \theta_3, \theta_4 \}$.   By Proposition~\ref{prop2.02} (b), $R$ has {\offdiag} (resp. $R$ has {\offdr}).  It follows from our hypothesis that $\{ \theta_1, \theta_2, \theta_3, \theta_4\}$ are either co-linear or co-circular.   Since this is true for an arbitrary collection of four elements from $\sigma(T)$, we conclude that all of the eigenvalues of $T$ are either co-linear or co-circular.  As before, this implies the existence of $\alpha, \beta \in \bbC$ and either a hermitian operator $L$ or a unitary operator $V$ such that $T = \alpha I + \beta L$, or $T = \alpha I + \beta V$.
\end{itemize}

We now concentrate on proving that a $4 \times 4$ matrix $T$ has {\offdiag} (resp. {\offdr}) if and only if the eigenvalues of $T$ are either co-linear or co-circular.
\end{rem}

\begin{lem} \label{lem3.12}
Let $X, Y \in \bbM_2(\bbC)$ and suppose that $\norm X \norm_2 = \norm Y \norm_2$.   The following are equivalent:
\begin{enumerate}
	\item[(a)]
	$\norm X \norm  = \norm Y \norm$; 
	\item[(b)]
	$\mathrm{tr}( (X^* X)^2) = \mathrm{tr}( (Y^* Y)^2)$ ;
	\item[(c)]
	$| \det (X)| = | \det (Y) |$.
\end{enumerate}
\end{lem}

\begin{pf}
Again, since the Fr\"{o}benius norm, the operator norm, and the trace functional are all invariant under unitary conjugation, we may assume without loss of generality that $X^* X$ and $Y^* Y$ are not only positive but diagonal, say 
\[
X^* X = \begin{bmatrix} x_1 & 0 \\ 0 & x_2 \end{bmatrix}, \ \ \ \ \ Y^* Y = \begin{bmatrix} y_1 & 0 \\ 0 & y_2 \end{bmatrix}, \]
with $0 \le x_1, x_2, y_1, y_2$.

The hypothesis that $\norm X \norm_2 = \norm Y \norm_2$ is the statement that $\varrho = x_1 + x_2 = y_1 + y_2$.  

\begin{enumerate}
	\item[(a)] implies (b).
	
	Suppose that $\norm X \norm  = \norm Y \norm$.   Then $\norm X \norm^2 = \norm Y \norm^2$ and so  $\max \{ x_1, x_2 \} = \max \{ y_1, y_2 \}$.  By reindexing if necessary, we may assume that $x_1 = y_1$.    But we have also assumed that $x_1 + x_2 = y_1 + y_2$, and so $x_2 = y_2$.   It follows that 
\[
\mathrm{tr} ( (X^* X)^2) = x_1^2 + x_2^2 = y_1^2 + y_2^2 = \mathrm{tr} ((Y^*Y)^{2}). \]
	\item[(b)] implies (c).
	
	Our current hypotheses are that $x_1 + x_2 = y_1 + y_2$ and that $x_1^2 + x_2^2 = y_1^2 + y_2^2$.   Thus 
	\begin{align*}
	|\det(Y)|^2 
		&= \det (Y^* Y)\\
		&= y_1 y_2 \\
		&= \frac{1}{2} \left((y_1 + y_2)^2 - (y_1^2 + y_2^2)\right)\\ 
		&= \frac{1}{2} \left((x_1 + x_2)^2 - (x_1^2 + x_2^2)\right)\\ 
		&= x_1 x_2\\ 
		&= \det (X^* X)\\ 
		&= |\det(X)|^2,
	\end{align*}	
from which (c) follows.	
	\item[(c)] implies (a).
	
	Suppose that $| \det(X)| = |\det (Y)|$.   Then, as we have just computed, $x_1 x_2 = | \det(X)|^2 = |\det (Y)|^2 = y_1 y_2$.
	
	But then $x_1 + x_2 = y_1 + y_2$ and $x_1 x_2 = y_1 y_2$ together imply that $\{ x_1, x_2 \} = \{ y_1, y_2 \}$.   In particular, 
	\[
	\norm X \norm^2 = \norm X^* X \norm = \max \{ x_1, x_2 \} = \max \{ y_1, y_2 \} = \norm Y^* Y \norm = \norm Y\norm^2. \]
	This completes the proof.
\end{enumerate}	
\end{pf}


\begin{thm}\label{T1}
Suppose that $T$ is an invertible normal block-matrix in $\bbM_{4}(\bbC)$ with $2\times 2$ blocks, and 
the off-diagonal corners of $T$ have equal rank (respectively equal operator norm), then the same is true for the off-diagonal corners of   $T^{-1}$.
\end{thm}

\begin{pf}
Let us start with the case of equal ranks, and employ a proof by contradiction, supposing that $\displaystyle T=\left[\begin{smallmatrix}
 A & B \\ C & D 
\end{smallmatrix}\right]$ is an invertible normal matrix with $\displaystyle T^{-1}=\left[\begin{smallmatrix}
 A' & B' \\ C' & D' 
\end{smallmatrix}\right]$, and $\mathrm{rank}\, {B}=\mathrm{rank}\, {C}, \text{  but  } \mathrm{rank}\, {B'}\neq\mathrm{rank}\, {C'}.$

Since $T^{-1}$ is normal, every invariant subspace of $T^{-1}$ is reducing, and so if either ${B'}$ or ${C'}$ is zero then both ${B'}$ and ${C'}$ are zero, contradicting our hypothesis. Hence we may assume that one of $B'$ and $C'$ has rank $1$ and the other has rank equal to $2$. 
Passing to $T^{{*}}$ if necessary, we can assume without loss of generality that 
$\mathrm{rank}\, {C'}=1<2=\mathrm{rank}\, {B'}.$

In particular, $B'$ is invertible, as is $T^{-1}$, and therefore 
$\displaystyle B=\left(T^{-1}|B'\right)^{-1}$,
so that $B$ is invertible.  Consequently $C$ has rank $2$ and is invertible. Hence 
$\displaystyle C'=\left(T|C\right)^{-1}$,
and therefore $C'$ is invertible, i.e. has rank $2$, equal to that of $B'$, contradicting our hypothesis.

Next let us deal with the case of equal operator norms. Let us suppose that 
$\Norm{B}=\Norm{C}$,
or equivalently, by Lemma \ref{lem3.12}, that
$|\det{B}|=|\det{C}|$.

First, let us treat the case ``$A$ is invertible''. In this case\[%
\det{C'}=\det{\left(-(T|A)^{-1}CA^{-1}\right)}=\frac{\det{C}}{\det{(T|A)}\det{A}}
\]%
and\[%
\det{B'}=\det{\left(-A^{-1}B\,(T|A)^{-1}\right)}=\frac{\det{B}}{\det{(T|A)}\det{A}},
\]%
so that $\det{C'}=\det{B'}$, and therefore
$\Norm{B'}=\Norm{C'}$,
again by Lemma \ref{lem3.12}.

Now, the remaining case is ``$A$ is not invertible''. In this case there is a sequence $\left[\alpha_{k}\right]_{k\in\N}$ convergent to zero and such that each $\alpha_{k}$ is neither an eigenvalue of $T$, nor of $A$. Applying the already settled case ``A is invertible'' to each (invertible and normal) $T-\alpha_{k}\, I$, {we can conclude that for each $k$ the off-diagonal corners of $\left(T-\alpha_{k}\,  I\right)^{-1}$
have equal norms. Yet $\lim_{k \to \infty} \norm (T-\alpha_k I)^{-1} - T^{-1} \norm = 0$, 
and therefore the off-diagonal corners of $\left(T-\alpha_{k}\,  I\right)^{-1}$ converge to those of $T^{-1}$, showing that the latter have equal norms as well, and the proof is complete.}
%
%
\end{pf}

\begin{cor}\label{C1}
Suppose that $T$ is a normal block-matrix in $\bbM_{4}(\bbC)$ with $2\times 2$ blocks such that a M\"{o}bius map $$M(z)=\frac{az+b}{cz+d}$$ is finite at all eigenvalues of $T$. 

If the off-diagonal corners of $T$ have equal rank (respectively, equal operator norm), then the same is true for 
the off-diagonal corners of $M(T)$.
\end{cor}
\smallskip

\begin{pf}
The claim is obviously true if $c=0\neq d$. Let us consider the case $c\neq 0.$ In this case, $\frac{-d}{c}$ is not an eigenvalue of $T$, and \[%
M(z)=\frac{a}{c}+\left(\frac{b-\frac{ad}{c}}{c}\right)\cdot \frac{1}{z+\frac{d}{c}}.
\]%
If the off-diagonal corners of $T$ have equal rank (respectively, equal operator norm), then the same is true for $T+\frac{d}{c}I$. Then, by Theorem \ref{T1}, the off-diagonal corners of $\left(T+\frac{d}{c}I\right)^{-1}$ have equal rank (respectively, equal operator norm), and so the same can be said about the  off-diagonal corners of $M(T)$.
\end{pf}

\begin{cor}\label{C2}
If $T\in\Mn{4}(\C)$ has {\offdiag} or if $T$ has {\offdr}, then $T$ is normal and $M(T)$ has the same property for any M\"{o}bius map $M$ that is finite on the spectrum of $T$.
\end{cor}\smallskip

\begin{pf}
This is the consequence of Proposition \ref{prop3.02}, Corollary \ref{C1} and the standard analytic functional-calculus fact that\[%
M\left(U^{^{*}}TU\right)=U^{^{*}}M(T)U.
\]%
\end{pf}

\begin{prop}\label{T2}
If $T$ is a normal block-matrix in $\bbM_{4}(\bbC)$ with $2\times 2$ blocks, and the spectrum of $T$ is $\{0,1,2, \delta\}$, where $\delta\notin\R$, then there exists a unitary block-matrix $U$ in $\bbM_{4}(\bbC)$ such that, with respect to the $2\times 2$ block partitioning,
\[%
\mathrm{rank}\, \left(U^{^{*}}TU\right)_{_{NE}}<2=\mathrm{rank}\, {\left(U^{^{*}}TU\right)_{_{SW}}}.
\]%
\end{prop}

\begin{pf} Every complex number $\delta$ other than 2 can be expressed as $2-8/(6+\beta)$ for a unique $\beta\neq-6$. Furthermore, $\delta$ is real exactly when $\beta$ is real.
 
Hence, after applying a unitary similarity we can assume without loss of generality that \[%
T=\left(
\begin{array}{cccc}
 0 & 0 & 0 & 0 \\
 0 & 1 & 0 & 0 \\
 0 & 0 & 2 & 0 \\
 0 & 0 & 0 & 2-\frac{8}{\beta+6}
   \\
\end{array}
\right)
\]%
The unitary $U$ shall be the product $VW$ of the unitaries \[%
V=\left(
\begin{array}{cccc}
 \sigma  & 0 & \gamma  & 0 \\
 0 & \gamma  & 0 & \sigma  \\
 \gamma  & 0 & -\sigma  & 0 \\
 0 & \sigma  & 0 & -\gamma  \\
\end{array}
\right)\ \AND\ W=\left(
\begin{array}{cccc}
 \frac{1}{\sqrt{2}} &
   \frac{e^{i \theta
   }}{\sqrt{2}} & 0 & 0 \\
 \frac{1}{\sqrt{2}} &
   -\frac{e^{i \theta
   }}{\sqrt{2}} & 0 & 0 \\
 0 & 0 & \frac{1}{\sqrt{2}} &
   \frac{1}{\sqrt{2}} \\
 0 & 0 & \frac{1}{\sqrt{2}} &
   -\frac{1}{\sqrt{2}} \\
\end{array}
\right),
\]%
where $\sigma, \gamma$ and $\theta$, to be specified later, are subject to the conditions:\[%
\sigma\neq \gamma,\ \ \sigma^{2}+\gamma^{2}=1,\ \ 0<\sigma,\ \ 0<\gamma, \AND \theta \text{ is not an integer multiple of } \pi.
\]%
A direct calculation shows that \[%
\left(U^{^{*}}TU\right)_{_{NE}}=\left(
\begin{array}{cc}
 -\frac{(3 \beta+10) \gamma 
   \sigma }{2 (\beta+6)} &
   -\frac{4 \gamma
   ^2}{\beta+6}-\frac{1}{2} e^{i
   \theta } \sigma ^2 \\
 -\frac{1}{2} e^{-i \theta }
   \gamma ^2-\frac{4 \sigma
   ^2}{\beta+6} & -\frac{(3 \beta+10)
   \gamma  \sigma }{2 (\beta+6)}
   \\
\end{array}
\right)
\]%
and\[%
\left(U^{^{*}}TU\right)_{_{SW}}=\left(
\begin{array}{cc}
 -\frac{(3 \beta+10) \gamma 
   \sigma }{2 (\beta+6)} &
   -\frac{1}{2} e^{i \theta }
   \gamma ^2-\frac{4 \sigma
   ^2}{\beta+6} \\
 -\frac{4 \gamma
   ^2}{\beta+6}-\frac{1}{2} e^{-i
   \theta } \sigma ^2 &
   -\frac{(3 \beta+10) \gamma 
   \sigma }{2 (\beta+6)} \\
\end{array}
\right),
\]%
and that
\begin{align*}
\det{\left(U^{^{*}}TU\right)_{_{NE}}}
&=\left(\frac{-2\gamma ^2\sigma ^2}{\beta+6}\right)
\left(\frac{1}{\frac{e^{i \theta }\sigma ^2}{\gamma ^2}} +\frac{e^{i \theta }\sigma ^2}{\gamma ^2}-\beta\right),
\end{align*}
while\[%
\det{\left(U^{^{*}}TU\right)_{_{SW}}}=\left(\frac{-2}{\beta+6}\right)\left(\gamma ^4 e^{i \theta
   }+e^{-i \theta } \sigma
   ^4-\gamma ^2 \sigma ^2
   \beta\right).
\]%
Now one can see that\[%
\det{\left(U^{^{*}}TU\right)_{_{NE}}}-\det{\left(U^{^{*}}TU\right)_{_{SW}}}=\frac{4 i \left(\gamma
   ^4-\sigma ^4\right) \sin
   (\theta )}{\beta+6}\neq 0,
\]%
because of the conditions that we have imposed on $\sigma, \gamma$ and $\theta$.

Note that $\displaystyle\frac{e^{i \theta }\sigma ^2}{\gamma ^2}$ can take on any non-real complex value even when $\sigma, \gamma$ are restricted to be distinct positive numbers whose squares add up to $1$, and $\theta$ is not an integer multiple of $\pi$.

It is also easy to see that the equation 
$\displaystyle\zeta+\frac{1}{\zeta}=\beta$ 
has a complex solution for $\zeta$, and since $\beta\not\in\R$, the solution cannot be real.

It follows that there exist $\sigma, \gamma$ and $\theta$ satisfying the conditions:\[%
\sigma\neq \gamma,\ \ \sigma^{2}+\gamma^{2}=1,\ \ 0<\sigma,\ \ 0<\gamma, \AND \theta \text{ is not an integer multiple of } \pi,
\]%
as well as the condition \[%
\frac{1}{\frac{e^{i \theta }\sigma ^2}{\gamma ^2}} +\frac{e^{i \theta }\sigma ^2}{\gamma ^2}=\beta.
\]%
These are the $\sigma, \gamma$ and $\theta$ that we use in the construction of $U$, and it is now clear that such a $U$ is the one we seek, since for this $U$:\[%
0=\det{\left(U^{^{*}}TU\right)_{_{NE}}}\neq \det{\left(U^{^{*}}TU\right)_{_{SW}}}.
\]%
\end{pf}


\begin{cor}\label{cor3.15}
If $T\in\Mn{4}(\C)$ has {\offdiag} or if $T$ has {\offdr}, then $T$ is normal and the spectrum of $T$ is circlinear.
\end{cor}\smallskip

\begin{pf}
Such $T$ has to be normal by Proposition \ref{prop3.02}. To verify the rest of the claims we proceed by contradiction. Suppose that the eigenvalues $\lambda_{0}, \lambda_{1}, \lambda_{2}, \lambda_{3}$ of $T$ are not circlinear. Then they are all distinct. 

By the spectral mapping theorem, given a M\"{o}bius map $M$ that is finite on the spectrum of $T$, the eigenvalues of $M(T)$ are the images of the eigenvalues of $T$ under $M$. It is well-known that M\"{o}bius maps take circlines to circlines, and exhibit sharp three-fold transitivity. 

In particular there is a unique M\"{o}bius map ${M}_{_{o}}$ such that\[%
{M}_{_{o}}(\lambda_{i})=i,\ \text{ for } i=0,1,2.
\]%
If ${M}_{_{o}}$ sends $\lambda_{3}$ to $z_{_{o}}$ that is a real number or ``$\infty$'', then the inverse of ${M}_{_{o}}$ sends $0,1,2, z_{_{o}}$ to 
$\lambda_{0}, \lambda_{1}, \lambda_{2}, \lambda_{3}$, indicating that the latter set is part of the image (under ${M}_{_{o}}$) of the extended real line, and hence must be circlinear, contrary to our hypothesis. Therefore ${M}_{_{o}}$ sends $\lambda_{3}$ to some complex non-real number $\delta$.

Applying Corollary \ref{C2} and Proposition \ref{T2} to ${M}_{_{o}}(T)$ yields a contradiction, and the proof is complete.
\end{pf}



By combining Remark~\ref{rem3.08} and Corollary~\ref{cor3.15}, we obtain the main theorem of this section.

\smallskip

\begin{thm} \label{thm3.16}
Let $n \ge 4$ be an integer and $T \in \bbM_n(\bbC)$.  The following are equivalent.
\begin{enumerate}
	\item[(a)]
	$T$ has {\offdiag}.
	\item[(b)]
	$T$ has {\offdr}.
	\item[(c)]
	One of the following holds.
	\begin{enumerate}
		\item[(i)]
		There exist $\lambda, \mu \in \bbC$ and $V \in \bbM_n(\bbC)$ unitary such that $T = \lambda I + \mu V$. 
		\item[(ii)]
		There exist $\lambda, \mu \in \bbC$ and $L = L^*\in \bbM_n(\bbC)$ such that 
		$T = \lambda I +\mu L$.
	\end{enumerate}
\end{enumerate}	
\end{thm}

\subsection{} \label{sec3.16}
{There is also an alternative proof for Theorem \ref{thm3.16} that does not involve M\"{o}bius maps, and while we have chosen not to include it here, we will gladly share it with an interested reader. Clearly the invariance of {\offdiag} and {\offdr} under M\"{o}bius maps (as in Corollary \ref{C2}) can  be inferred from Theorem \ref{thm3.16}. }


\section{The infinite-dimensional setting -- {\offdiag}} \label{sec4}



\subsection{} \label{sec4.01}
Let us now consider the case where the underlying Hilbert space is infinite-dimensional and separable.  We begin by studying operators with {\offdiag}.   In the finite-dimensional setting, we saw that any such operator is normal.  While this is also true in the infinite-dimensional setting, the proof is rather different.


\bigskip 

Recall that an operator $T \in \bofh$ is said to be \textbf{strongly reductive} if, whenever $(P_n)_{n=1}^\infty$ is a sequence of orthogonal projections such that $\lim_n \norm P_n^\perp T P_n \norm = 0$, it follows that $\lim_n \norm P_n T P_n^\perp \norm = 0$ (or equivalently, $\lim_n \norm P_n T - T P_n \norm = 0$).  Let us say that a compact set $\Omega \subseteq \bbC$ is \textbf{Lavrentiev} if it has empty interior and if $\bbC \setminus \Omega$ is connected.

\begin{prop} \label{prop4.02}
Let  $\hilb$ be an infinite-dimensional, separable Hilbert space.  If $T \in \bofh$ has {\offdiag}, then $T$ is normal and has Lavrentiev spectrum.
\end{prop}

\begin{pf}
It is an immediate consequence of the definition that if $T$ has {\offdiag}, then $T$ is strongly reductive.   
It was shown by Harrison~\cite{Har1975} that any strongly reductive operator has Lavrentiev spectrum.   Apostol, Foia\c{s} and Voiculescu~\cite{AFV1976} showed that any strongly reductive operator is normal.

It is also easy to obtain the normality of $T$ which enjoys {\offdiag} directly.
Suppose that $T \in \bofh$ has {\offdiag}, and let $e \in \hilb$ be an arbitrary vector of norm one.  Then $P_e (x) = \langle x, e \rangle e$, $x \in \hilb$ defines a rank one projection.   By hypothesis, $\norm P_e^\perp T^* P_e \norm = \norm P_e T P_e^\perp \norm = \norm P_e^\perp T P_e \norm$.

Now 
\begin{align*}
\langle T T^* e, e \rangle = \norm T^* e \norm^2 
	&= \norm P_e T^* e \norm^2 + \norm P_e^\perp T^* e \norm^2 \\
	&= \norm P_e T^* P_e e \norm^2 + \norm P_e^\perp T^* P_e e \norm^2 \\
	&= \norm P_e T^* P_e \norm^2 + \norm P_e^\perp T^* P_e \norm^2.
\end{align*}	

Similarly, 
\[
\langle T^* T e, e \rangle = \norm P_e T P_e\norm^2 + \norm P_e^\perp T P_e\norm^2. \]

Now 
\[
\norm P_e T P_e \norm  = \norm P_e T^* P_e \norm, \]
and combining this with the fact that $\norm P_e^\perp T^* P_e \norm  = \norm P_e^\perp T P_e \norm$ from above, we see that 
\[
\langle T T^* e, e \rangle = \langle T^* T e, e \rangle. \]
Since $e$ was an arbitrary norm-one vector, we conclude that $T T^* = T^* T$; i.e. that $T$ is normal.

\end{pf}

%
%


\subsection{} \label{sec4.03}
In Section~\ref{sec2}, we noted that if $L = L^* \in \bofh$ and $\lambda, \mu \in \bbC$, then $\lambda I + \mu L$ has {\offdiag}.   Although in the finite-dimensional setting every unitary operator $V$ also has {\offdiag}, this is no longer true in the infinite-dimensional setting, as the following counterexample shows.

	Let  $\{ e_n \}_{n \in \bbZ}$ be an orthonormal basis for our Hilbert space $\hilb$, and consider the \textbf{bilateral shift} operator $W$ determined by $W e_n = e_{n-1}$, $n \in \bbZ$.   Then $W$ is unitary and $\sigma(W) = \bbT$.  By Proposition~\ref{prop4.02}, $W$ does not have {\offdiag}.   

\smallskip
	
	This can be seen directly as well.  If $P$ is the orthogonal projection of $\hilb$ onto $\cM = \ol{\mathrm{span}} \{ e_n: n \le 0\}$, then $\cM$ is invariant for $W$, so that $P^\perp W P = 0$.   However, $0 \ne e_0 = W e_1 = P W P^\perp e_1$, so that $P W P^\perp \ne 0$.  \emph{A fortiori}, $W$ has neither {\offdiag} nor {\offdr}. Furthermore, since the set of operators having {\offdiag} is clearly (norm-)closed, no operator close enough to $W$ has {\offdiag}.
	
	
	For $n \ge 3$, let $C_n \in \bbM_n(\bbC)$ denote an $n$-cycle; that is, there is an orthonormal basis $\{ e_k\}_{k=1}^n$ of $\bbC^n$, such that $C_n e_k = e_{k+1}$, $1 \le k \le n-1$, {and} $C_n e_n = e_1$.  It follows easily from the results in~\cite{Dav1986} that there exists a sequence $(V_n)_{n=1}^\infty$ of unitary operators with $V_n \simeq C_n \otimes I$ for all $n \ge 1$ such that $\lim_n V_n = W$.  {Thus $V_n$ does not have {\offdiag}
for all sufficiently large $n$.}   That is, for sufficiently large $n$, $C_n \otimes I$ fails to have {\offdiag}, despite the fact that $C_n \in \bbM_n(\bbC)$ has {\offdiag} by Proposition~\ref{prop3.03}.  In fact, as we shall soon see, $C_n \otimes I$ fails to have {\offdiag} for all $n \ge 3$.


\subsection{} \label{sec4.04}
\mbox{} {L}et us recall that the \textbf{numerical range} of $T \in \bofh$ is the set $W(T) = \{ \langle T x, x \rangle: \norm x \norm = 1 \}$, and the \textbf{numerical radius} of $T$ is $w(T) = \sup \{ |\lambda|: \lambda \in W(T) \}$.  It is known that the numerical range of an operator is always a convex set (this is the classical Toeplitz-Hausdorff Theorem), and that the closure of the numerical range of $T$ always contains $\sigma(T)$ { (see, for example, Problem~214 of~\cite{Hal1982})}. 

{ In trying to characterize operators with {\offdiag}, we  know that we may restrict our attention to normal operators. For a normal }$M$, it is known that $\ol{W(M)} = \mathrm{co} (\sigma(M))$, that is, the closure of the numerical range of $M$ is the convex hull of the spectrum of $M$. { If $\sigma(M)$ happens to be finite, then (by the Spectral Theorem) $\sigma(M)$ consists of the eigenvalues of $M$, and these belong to the numerical range of $M$, as does their convex hull. Hence in such a case $W(M)=\mathrm{co} (\sigma(M))$.
}

{ The} numerical radius defines a norm on $\bofh$ which is equivalent to the operator norm, { because} $\frac{1}{2} \norm T \norm \le w(T) \le \norm T \norm$ for all $T \in \bofh$.   (See, e.g.~\cite{Hal1982}, Chapter 22 for all of these results.)

{ A} \textbf{state} on $\bofh/\kofh$ is a positive linear functional of norm one.  For $T \in \bofh$, the \textbf{essential numerical range} { $W_e(T)$} of $T$ is the set $\{ \varphi(\pi(T)): \varphi \mbox{ is a state on } \bofh/\kofh\}$.   It is known that { $W_e(T)$ is} closed and convex, and so it follows that  $W_e(T)=\mathrm{co} (\sigma_{e}(T))$, whenever $T$ is normal. 

\begin{thm} \emph{(}\textbf{Fillmore-Stampfli-Williams}; Theorem 5.1 of~\cite{FSW1972}\emph{)}  \label{thm4.05}
For $T \in \bofh$, the following conditions are equivalent:
\begin{enumerate}
	\item[(a)]
	$0 \in W_e(T)$.
	\item[(b)]
	There exists an orthonormal sequence $(e_n)_{n=1}^\infty$ in $\hilb$ {such} that $\lim_n \langle T e_n, e_n \rangle = 0$.
	\item[(c)]
	$0 \in \cap \{ \ol{W(T+F)} : F \mbox{ is of finite rank}\}$.
\end{enumerate}	
\end{thm}

From this it easily follows that $W_e(T) = \cap \{ \ol{W(T+F)} : F \mbox{ is of finite rank}\}$.


\subsection{} \label{sec4.06}  
 If $R \in \cB(\cK)$ where $\cK \subseteq \cH$ is a subspace of $\cH$, then $T$ is said to be a \textbf{dilation} of $R$ if, relative to the decomposition $\cH = \cK \oplus \cK^\perp$, we may write 
\[
T = \begin{bmatrix} R & B \\ C & D \end{bmatrix} \]
for some choice of $B, C$ and $D$.

{ Recall} that if $\{ e_n \}_{n=1}^\infty$ is an orthonormal basis for $\hilb$, then the \textbf{unilateral forward shift} on $\hilb$ is the operator $S \in \bofh$ satisfying $S e_n = e_{n+1}$ for all $n \ge 1$.

{   The following result of Choi and Li will be useful.  }


\begin{thm} \emph{(}\textbf{Choi-Li}; Theorem 4.3 of~\cite{CL2001}\emph{)} \label{thm4.07}
Suppose that $A \in \bofh$, and $T \in \bbM_3(\bbC)$ has a non-trivial reducing subspace.   {Then  $A$ has a dilation that is unitarily equivalent to $T \otimes I$ if and only if $W(A) \subseteq W(T)$.}
\end{thm}


\begin{thm} \label{thm4.08}
Suppose that $V \in \bbM_3(\bbC)$ is a unitary operator and that $0$ lies in the interior of $W(V)$.   Then $V \otimes I$ does not have {\offdiag}.
\end{thm}

\begin{pf}
Clearly every unitary operator $V$ in $\bbM_3(\bbC)$ has a non-trivial reducing subspace.  Note also that if $S \in \bofh$ is the unilateral forward shift, $\mathrm{spr}(S) \le \norm S \norm = 1$.  (In fact, $\mathrm{spr}(S) = 1$, but that is not important here.)

Thus $\mathrm{spr} (\eps S) \le \eps$ for all $\eps > 0$.   Since $0$ lies in the \emph{interior} of $W(V)$, there exists $\eps_0 > 0$ such that $W(\eps_0 S) \subseteq W(V)$.  Let $A = \eps_0 S$.   By Theorem~\ref{thm4.07} above, we may write 
\[
V \otimes I \simeq \begin{bmatrix} A & B \\ C & D \end{bmatrix}. \]

Let $P = S S^*$, so that $P$ is an orthogonal projection of co-rank one (i.e. the rank of $(I-P)$ is equal to one).   In particular, $\sigma(P) = \{ 0, 1 \}$.   

\smallskip

Then, since $V \otimes I$ is unitary, we have that 
\begin{align*}
B B^* = I - A A^* = I - \eps_0^2 P \\
C^* C = I - A^* A = I - \eps_0^2 I.
\end{align*}

It follows that 
\[
\norm B \norm^2  = \norm B B^* \norm = \norm I - \eps_0^2 P \norm = 1, \]
while 
\[
\norm C \norm^2 = \norm C^* C \norm = \norm I - \eps_0^2 I \norm  = 1-\eps_0^2. \]

In particular, $\norm B \norm \ne \norm C \norm$, so that $V \otimes I$ does not have {\offdiag}.
\end{pf}

It follows from Theorem \ref{thm4.08} that $C_n \otimes I$ does not have {\offdiag} for $n \ge 3$, since in such a case $0$ lies in the interior of $W(C_n)$.


\smallskip

\begin{cor} \label{cor4.09}
Suppose that $U \in \bofh$ is a unitary operator  that has {\offdiag}. Then $0$ does not lie in the interior of $W_e(U)$.
\end{cor}

\begin{pf} We prove a contrapositive implication. Suppose that $0$ lies in the interior of the essential numerical range of $U$, that is in the interior of the convex hull of the essential spectrum of $U$.
Then there exist $\alpha, \beta$ and $\gamma$ in the essential spectrum of $U$ such that $0$ lies in the interior of the convex hull of $\{ \alpha, \beta, \gamma\}$.
Let $V \in \bbM_3(\bbC)$ be a unitary operator with spectrum $\{ \alpha, \beta, \gamma\}$.   Then, as noted in Section~\ref{sec4.04}, $W(V)$ is closed and $W(V) = \mathrm{co} \{ \alpha, \beta, \gamma\}$. Hence $0$ lies in the interior of $W(V)$.   By Theorem \ref{thm4.08}, $V \otimes I$ does not have {\offdiag}.

Since $\alpha, \beta$ and $\gamma$ lie in the essential spectrum of the normal operator $U$, $U$ is approximately unitarily equivalent to $U \oplus (V \otimes I)$.    (This is a consequence of the Weyl-von Neumann-Berg Theorem for normal operators -- see, e.g. Theorem~II.4.4 of~\cite{Dav1996} -- and can also be deduced from the results of~\cite{Dav1986}.)   It now follows from Proposition~\ref{prop2.03} (b) that $U$ does not have {\offdiag}.
\end{pf}


\begin{thm} \label{thm4.10}
Suppose that $U \in \bofh$ is a unitary operator and $0$ does not lie in $W_e(U)$.   Then $U$ has {\offdiag}.
\end{thm}

\begin{pf} We prove a contrapositive implication.
Let $U$ be a unitary operator which fails to have {\offdiag}.   Then there exists a projection $P \in \bofh$ such that with respect to the decomposition $\hilb = P\hilb \oplus P^\perp \hilb$, we may write 
\[
U = \begin{bmatrix} A & B \\ C & D \end{bmatrix} \]
where $\norm B \norm \ne \norm C \norm$.   

As noted in Remark~\ref{rem2.04}, this can only happen if one of the following holds:   
\begin{enumerate}
	\item[(a)]
	either $0 \in \sigma(A A^*) $ but $0 \not \in \sigma (A^* A)$, 
	\item[(b)]
	$0 \in \sigma (A^* A)$ but $0 \not \in \sigma(A A^*)$.
\end{enumerate}

Since an operator $T$ has {\offdiag} if and only if $T^*$ has {\offdiag}, by replacing $U$ by $U^*$ if necessary (which does not affect the conclusion, as $0 \in W_e(U)$ if and only if $0 \in W_e(U^*)$), we may assume without loss of generality that $0 \in \sigma(A A^*)$ but $0 \not \in \sigma(A^* A)$. In particular, $A$ is not invertible.  

Since $A^* A$ is invertible, we see that $A$ is bounded below.  Hence $\mathrm{ran}\, A$ is closed and $\mathrm{nul}\, A = 0$.  Therefore $A$ is semi-Fredholm.   	If $\mathrm{nul}\, A^* = 0$, then $\mathrm{ran}\, A = P\hilb$, so that $A$ is invertible, which is a contradiction.    Thus $\mathrm{nul}\, A^* > 0$, and so $\mathrm{ind}\, A < 0$.     

Since $\mathrm{ind}\, (A+F) = \mathrm{ind}\, A < 0$ for all finite-rank operators $F \in \cB(P \hilb)$, and since $\sigma(T) \subseteq \ol{W(T)}$ for all operators $T$, we may apply Theorem~\ref{thm4.05} of Fillmore, Stampfli and Williams to obtain: 
\begin{align*}
0 
	&\in \cap \{ \sigma(A+F) : F \in \cB(P \hilb), \ F \mbox{ finite-rank}\} \\
	&\subseteq \cap \{ \ol{W(A+F)}: F \in \cB(P \hilb),\  F \mbox{ finite-rank} \} \\
	&= W_e(A), 
\end{align*}	
Since $W_e(A) \subseteq W_e(U)$, the result follows.
\end{pf}


\begin{thm} \label{thm4.11}
Suppose that $U \in \bofh$ is unitary and that $0$ lies on the boundary of $W_e(U)$.  Then $U$ has {\offdiag}.   
\end{thm}

\begin{pf}
The hypotheses of the theorem imply that $\sigma_e(U)$ lies on a closed half-circle of $\bbT$, and  includes two diametrically opposite points.   By multiplying $U$ by an appropriate $\mu \in \bbT$ (which does not affect the conclusion of the Theorem), we may assume without loss of generality that $\sigma_e(U) \subseteq \bbT \cap \{ z\in \bbC: \mathrm{Re}{z} \ge 0\}$, and that $\{ i, -i \} \subseteq \sigma_e(U)$.

For each $n \ge 1$, let $\cC_n = \{ z \in \bbT: z= e^{i \theta}, \frac{\pi}{2} - \frac{1}{n} \le \theta \le \frac{\pi}{2} + \frac{1}{n}\}$, and let $Q_n$ be the spectral projection for $U$ corresponding to $\cC_n$.  Then $U = X_n + Y_n$, where $X_n = U Q_n^\perp$,  and where $Y_n = U Q_n$ is a unitary with $\sigma(Y_n) \subseteq \cC_n$.   Set $V_n = X_n + e^{i (\frac{\pi}{2}-\frac{1}{n})} Q_n$.

It is reasonably straightforward to check that $\norm U - V_n \norm = \norm Y_n - e^{i (\frac{\pi}{2}-\frac{1}{n})} Q_n \norm \le \frac{4 \pi}{n}$ and that $\sigma_e(V_n) \subseteq \Omega_n = \{ z= e^{i \theta} \in \bbT: -\frac{\pi}{2} \le \theta \le \frac{\pi}{2} - \frac{1}{n} \}$.

   Thus $0$ does not lie in the closed, convex hull of $\sigma_e(V_n)$, and in particular, $0$ does not lie in $W_{e}(V_n)$ for any $n \ge 1$.  By Theorem~\ref{thm4.10}, $V_n$ has {\offdiag}.  But as we saw in Proposition~\ref{prop2.03}, the set $\mathfrak{G_{norm}}$ of operators with {\offdiag} is closed, and thus $U$ has {\offdiag}.
\end{pf}

%
%
%
%
 	

Combining these results, and keeping in mind that $W_e(T)=\mathrm{co} (\sigma_{e}(T))$ for normal $T$, we obtain the following. 

\begin{cor} \label{cor4.12}
The following are equivalent for a unitary operator $U\in \bofh$.
\begin{enumerate}
	\item[(a)]
	$U$ has {\offdiag}.
	\item[(b)] $0$ does not lie in the interior of $W_e(U)$.
	\item[(c)] There exists a half-circle $\cC$ of $\bbT$ such that $\sigma_e(U) \subseteq \cC$\newline \textup{(}i.e.
	there exists $\mu \in \bbT$ such that $\sigma_e(U) \subseteq \bbT \cap \{ \mu z \in \bbC: \mathrm{Re}(z) \ge 0\})$. 
\end{enumerate}
\end{cor} 

%


We are now ready to state and prove the main theorem of this section.

\begin{thm} \label{thm4.13}
Let $T \in \bofh$.  The following conditions are equivalent.
\begin{enumerate}
	\item[(a)]
	$T$ has {\offdiag}.
	\item[(b)]
	One of the following holds. 
	\begin{enumerate}	
		\item[(i)]
		There exist $\lambda, \mu \in \bbC$ and $L=L^* \in \bofh$ such that $T = \lambda I + \mu L$.  
		\item[(ii)]
		There exist $\lambda, \mu \in \bbC$ with $\mu \ne 0$ and a unitary operator $U \in \bofh$ with $\sigma_e(U) \subseteq \bbT \cap \{ z \in \bbC: \mathrm{Re} (z) \ge 0\}$ such that $T = \lambda I + \mu U$.
	\end{enumerate}
\end{enumerate}
\end{thm}

\begin{pf}
Suppose first that (a) holds, and recall that this implies that $T$ is normal.

If $\sigma(T)$ has at most three points, then those points are either co-linear or co-circular.    It is routine to check from this that $T$ is either of the form of (b) (i), or there exist $\lambda, \mu \in \bbC$ with $\mu \ne 0$ and a unitary $U\in \bofh$ such that $T= \lambda I +\mu U$.  In this case, since $\mu \ne 0$, $T$ has {\offdiag} if and only if $U$ has {\offdiag}.  But then $U$ has {\offdiag} by our hypothesis on $T$, and so the spectral conditions on $U$ follow from Corollary~\ref{cor4.12}.

Thus we suppose that $\sigma(T)$ has cardinality at least 4, and we let $\{ \alpha, \beta, \gamma, \delta\}$ be four distinct points in $\sigma(T)$.   By Theorem II.4.4 in~\cite{Dav1996},  $T$ is approximately unitarily equivalent to an operator of the form $A \oplus D$, where $D = \mathrm{diag}  (\alpha, \beta, \gamma, \delta) \in \bbM_4(\bbC)$.  By Proposition~\ref{prop2.03} (b), $D$ has {\offdiag}.   By Corollary~\ref{cor3.15}, the eigenvalues of $D$ are either co-linear or co-circular.   Since this is true for any choice of four distinct points of $\sigma(T)$, we see that $\sigma(T)$ is either contained in a line -- in which case it is easily seen that there exist $\lambda, \mu$ and $L$ as in (b) (i) such that $T = \lambda I + \mu L$, or $\sigma(T)$ lies on a proper circle, i.e. there exist $\lambda, \mu \in \bbC$ with $\mu \ne 0$ and a unitary $U\in \bofh$ such that $T= \lambda I +\mu U$.  We argue as in the previous paragraph to obtain the spectral conditions on $U$.

\smallskip

Suppose next that (b) holds.  

If there exist $\lambda, \mu \in \bbC$ and $L=L^* \in \bofh$ such that $T = \lambda I + \mu L$, then $T$ has {\offdiag} by Proposition~\ref{prop2.02} (d).

If there exist $\lambda, \mu \in \bbC$ with $\mu \ne 0$ and a unitary operator $U \in \bofh$ with $\sigma_e(U) \subseteq \bbT \cap \{ z \in \bbC: \mathrm{Re} (z) \ge 0\}$ such that $T = \lambda I + \mu U$.  Then $T$ has {\offdiag} if and only if $U$ has {\offdiag} by Proposition~\ref{prop2.02} (a).  But $U$ has {\offdiag} by Corollary~\ref{cor4.12}, whence $T$ has {\offdiag}.
\end{pf}

\subsection{} \label{sec4.14}
\mbox{} {It} is natural to consider a weakening of the {\offdiag} obtained by restricting our attention to the finite-rank projections $P\in\bofh$ in the case when $\hilb$  is infinite-dimensional. As the reader can easily check, the results and proofs presented in this section readily demonstrate that such a ``weakening'' of the {\offdiag} is in fact equivalent to the original {\offdiag}.


\section{The infinite-dimensional setting -- {\offdr}} \label{sec5}



\subsection{} \label{sec5.01} 
We next turn our attention to the study of operators  with {\offdr}, acting on an infinite-dimensional Hilbert space.   Although we have not been able to obtain a complete classification of such operators, we will mention a number of interesting facts.


\bigskip 

We recall from above that an operator $X \in \bofh$ is said to be  (orthogonally) reductive  if for each projection $P \in \bofh$, the condition $P T P^\perp = 0$ implies that $P^\perp T P = 0$.     It is clear that if $X$ has {\offdr}, then $X$ must be reductive.   
			

It should be noted that not every normal operator is reductive.  Sarason~\cite{Sar1966} has shown that a normal operator $N$ is reductive if and only if $N^*$ lies in the weak operator topology closure of the set of polynomials in $N$.   As a concrete example, let $W \in \bofh$ be the bilateral shift; i.e. let $\{ e_n \}_{n\in \bbZ}$ be an orthonormal basis for $\hilb$ and let $W$ be defined by $W e_n = e_{n-1}$ for all $n \in \bbZ$.   It is well-known that $W$ is unitary with $\sigma(W) = \bbT$.   If $P$ is the orthogonal projection onto $\cM = \ol{\mathrm{span}} \{ e_n : n \le 0\}$, then clearly $\cM$ is invariant for $W$ -- so $\mathrm{rank}\, P^\perp W P = 0$, but it is easily verified that $P W P^\perp$ has rank $1$.   Thus $W$ fails to be reductive.   

The condition that an operator have {\offdr} is  strictly stronger than asking that it be orthogonally reductive -- see Example~\ref{eg5.09E} below.  It is worth observing that there is one inherent weakness in the definition of orthogonally reductive operators:  it is entirely possible that there might exist an operator with no non-trivial closed, invariant subspace, in which case the operator is reductive for trivial reasons.   On the other hand, it was shown by Popov and Tcaciuc~\cite{PT2013} that given any operator $T$ acting on an infinite-dimensional, complex, separable Hilbert space $\hilb$, there exists an orthogonal projection $P$ of infinite rank and co-rank such that $\mathrm{rank}\, P T P^\perp \le 1$.   (Their result actually holds for operators acting on reflexive Banach spaces and beyond, but we do not require that here.)  As such, {\offdr} always has significance for Hilbert space operators.

\bigskip


We begin with some observations regarding the general class ${\mathfrak{G}_{\mathfrak{rank}}}$ of operators with {\offdr}.

\begin{prop} \label{prop5.02}
Suppose that $T \in \bofh$ has {\offdr}.   Then $T$ is biquasitriangular; that is, $\mathrm{ind}\, (T-\lambda I) = 0$ for all $\lambda \in \rho_{sF}(T)$.
\end{prop}

\begin{pf}
It is clear that if $\lambda \in \bbC$, then $T - \lambda I$ and $(T-\lambda I)^*$ also have {\offdr}.   

Suppose that $\lambda \in \rho_{sF}(T)$ and that $\mathrm{ind}(T-\lambda I) \ne 0$.  By considering $(T-\lambda I)^*$ if necessary, we may assume that $\mathrm{ind}\, (T-\lambda I) > 0$ (it is possibly infinite).  

Thus $\nul\, (T-\lambda I) > \nul\, (T-\lambda I)^*$.  Write $\hilb = \ker\, (T-\lambda I) \oplus (\ker\, (T-\lambda I))^\perp$, and write 
\[
(T-\lambda I)  = \begin{bmatrix} 0 & B \\ 0 & D \end{bmatrix} \]
relative to this decomposition.  If $B = 0$, then $(T-\lambda I)^* = \begin{bmatrix} 0 & 0 \\ 0 & D^* \end{bmatrix}$, showing that $\nul\, (T-\lambda I)^* \ge \nul\, (T-\lambda I)$, a contradiction.  But then $B \ne 0$ implies that $T-\lambda I$ does not have {\offdr}, and hence neither does $T$.

The contrapositive is the statement that if $T$ has {\offdr}, then $T$ is biquastriangular.
\end{pf}


\subsection{} \label{sec5.03}
It is clear that if $T$ has {\offdr} and $\lambda$ is an eigenvalue for $T$, then it is a reducing eigenvalue for $T$; that is, we may write $T \simeq \lambda Q \oplus T_0$, where $Q$ is an orthogonal projection and $\lambda$ is no longer an eigenvalue for $T_0$ (though it may be an approximate eigenvalue for $T_0$).

By repeating this for each of the eigenvalues of $T$, this allows us to write $T \simeq M \oplus T_4$, where $M$ is a diagonal operator whose eigenvalues are precisely the eigenvalues of $T$, and where $T_4$ has no eigenvalues.   

	Since direct summands of operators with {\offdr} still have {\offdr}, it follows from the results of Section~\ref{sec3} that all of the spectrum of $M$ is either co-linear or co-circular.  
	That is, \emph{the eigenvalues of any operator} $T$ \emph{ with }{\offdr} \emph{ are either co-linear or co-circular, and they are reducing eigenvalues for } $T$.

	
Our next goal is to prove that every operator which satisfies {\offdr} is normal with circlinear spectrum.   We shall accomplish this through a sequence of lemmas.  It is worth noting that we shall not invoke the full strength of the {\offdr} hypothesis.  Indeed, for the next few results, we only require a weaker form of {\offdr} that requires that $T$ be reductive and that if $P$ is a projection for which $\mathrm{rank}\, P^\perp T P = 1$, then $\mathrm{rank}\, P T P^\perp = 1$.  It can in fact be shown that Proposition~\ref{prop5.02} also holds under this weaker hypothesis, though we shall not need that here.


\begin{prop} \label{prop5.04E}
Let $T \in \bofh$ and suppose that $T$ has {\offdr}.   Then there exist $\alpha, \beta, \gamma$ and $\delta \in \bbC$, not all equal to zero, and an operator $F \in \bofh$ of rank at most three such that 
\[
\alpha I + \beta T + \gamma T^* + \delta T^* T + F = 0. \]
\end{prop}

\begin{pf}
Fix $0 \ne \xi \in \hilb$.   We first claim that the set $S_\xi = \{ \xi, T \xi, T^* \xi, T^* T \xi \}$ is linearly dependent.

Let $\cM_\xi = \mathrm{span}\, \{ \xi, T \xi \}$.   If $\dim\, \cM_\xi = 1$, then clearly $\{ \xi, T \xi \}$ is linearly dependent, whence $S_\xi$ is linearly dependent and we are done.

Suppose therefore that $\dim\, \cM_\xi = 2$ and let $P_\xi$ denote the orthogonal projection of $\hilb$ onto $\cM_\xi$.    Note that $T \xi \in \cM_\xi$ implies that $\mathrm{rank}\, P_\xi^\perp T P_\xi \in \{0, 1\}$.   From our hypothesis, $\mathrm{rank}\, P_\xi^\perp T^* P_\xi \in \{0, 1\}$.

But then 
\begin{align*}
	\dim\, (P_\xi \hilb + P_\xi^\perp T^* P_\xi \hilb)
		&\le \dim\, (P_\xi \hilb)  +  \dim\, (P_\xi^\perp T^*  P_\xi \hilb) \\
		&\le 2  + 1 =3.
\end{align*}
Since $S_\xi \subseteq  P_\xi \hilb +P_\xi^\perp T^* P_\xi \hilb$, our claim follows.

As $0 \ne \xi \in \hilb$ was arbitrary, we see that the set $\{ I, T, T^*, T^* T \}$ is \emph{locally linearly dependent} in the sense of~\cite{Ami1965, Aup1988} and~\cite{MesS2002}.  By Theorem~2 of~\cite{Aup1988}, there exist $\alpha, \beta, \gamma,$ and $\delta \in \bbC$, not all equal to zero, such that 
\[
\mathrm{rank}\, (\alpha I + \beta T + \gamma T^* + \delta T^* T) \le 3. \]
This clearly implies the statement of the theorem.
\end{pf} 


We begin by dealing with the case where $\delta$ above is equal to zero.

\smallskip

\begin{lem} \label{lem5.05E}
Let $\hilb$ be a complex Hilbert space and suppose that $T \in \bofh$.   If there exist complex numbers $\alpha, \beta$ and $\gamma$, not all equal to zero, and $F \in \bofh$ of rank at most $m < \frac{1}{2} \, \dim \, \hilb$ such that 
\[
\alpha I + \beta T + \gamma T^* + F = 0, \]
then there exist a hermitian operator $R$, a finite-rank operator $L$ of rank at most $2 m$, and $\mu, \lambda \in \bbC$ such that 
\[
T = \lambda (R + L) + \mu I. \]
\end{lem}

\begin{pf}
\noindent{\textsc{Case 1.}} \ \ \ Suppose that $\gamma = 0$.

In this case, we have that $\alpha I  + \beta T + F = 0$.   If $\beta = 0$, then the fact that $F$ has finite rank $\mathrm{rank}\, F  = m <  \dim\, \hilb  = \mathrm{rank}\, I$ implies that $\alpha = 0 (= \beta = \gamma)$, contradicting our hypothesis.   Hence $\beta \ne 0$.

But then 
\[
T = - \alpha \beta^{-1} I - \beta^{-1} F. \]
\begin{itemize}
	\item{}
	If $\alpha = 0$, then 
	\[
	T = - \beta^{-1} F = -\beta^{-1} (0 + F) + 0 I \]
	expresses $T$ in the desired form.
	\item{}
	If $\alpha \ne 0$, then writing 
	\[
	T = - \alpha \beta^{-1} (I - \alpha^{-1} F) + 0 I \]
	expresses $T$ in the desired form.
\end{itemize}	
\bigskip
\noindent{\textsc{Case 2.}} \ \ \ Suppose that $\beta = 0$.

Then $\alpha I + \gamma T^* + F = 0$, and arguing as before, $\gamma \ne 0$.  Thus $T ^* = - \alpha \gamma^{-1} I - \gamma^{-1} F$.

\begin{itemize}
	\item{}
	If $\alpha = 0$, then 
	\[
	T^* =  -\gamma^{-1} (0 + F) + 0 I \]
	means that 
	\[
	T = -\ol{\gamma}^{-1} (0 + F^*) + 0 I \]
	expresses $T$ in the desired form.
	\item{}
	If $\alpha \ne 0$, then writing 
	\[
	T^* = - \alpha \gamma^{-1} (I - \alpha^{-1} F) + 0 I \]
	means that 
	\[
	T = -\ol{\alpha} \ol{\gamma}^{-1} (I - \ol{\alpha}^{-1} F^*) + 0 I \]
	expresses $T$ in the desired form.
\end{itemize}	

\bigskip
\noindent{\textsc{Case 3.}} \ \ \ Suppose that $\beta \ne 0 \ne \gamma$.

We have that $\alpha I + \beta T + \gamma T^* + F = 0$, whence $\ol{\alpha} I + \ol{\beta} T^* + \ol{\gamma} T + F^* = 0$.  Set $\rho = (\alpha + \ol{\alpha})$, $\theta = \beta + \ol{\gamma}$ and $F_0 = F + F^*$.   Adding the two previous equations involving $T$ yields:
\[
\rho I + \theta T + \ol{\theta} T^* + F_0 = 0, \]
and $\mathrm{rank}\, F_0 \le 2 \ \mathrm{rank}\, F \le 2 m < \dim\, \hilb.$

\bigskip
\noindent{\textsc{Subcase 3.A.}}\ \ \ $\theta = 0$.

If $\theta = 0$, then $\rho I + F_0 = 0$, combined with the fact that $\dim\, \hilb > \mathrm{rank}\, F_0$ implies that $\rho = 0 = F_0$.   That is, $\alpha \in i \bbR$ and $\gamma = -\ol{\beta} \ne 0$, so that 
\[
\alpha I + \beta T - \ol{\beta} T^* + F = 0. \]
Let $A  = \beta T + \frac{\alpha}{2} I$, and let $A = R+ i B$ be the Cartesian decomposition of $A$, so that $R$ and $B$ are hermitian.   The above equation shows that $0 = (A - A^*) + F = 2 i B + F$, and thus $B$ has finite rank at most $m$ and
\[
T = \beta^{-1} (R + i B) - \frac{\alpha \beta^{-1}}{2} I \]
expresses $T$ in the desired form.

\bigskip
\noindent{\textsc{Subcase 3.B.}}\ \ \ $\theta \ne 0$.

We have 
\[
\rho I + \theta T + \ol{\theta} T^* + F_0 = 0, \]
where $\mathrm{rank}\, F_0 \le 2 m$ and $\rho \in \bbR$.

Let $\dss \kappa = \frac{\rho i}{2 | \theta|^2} \in i \bbR$, and $A = (\ol{\theta} i)^{-1} T - \kappa I$.   Then $T = \ol{\theta} i(A + \kappa I)$ and our equation $\rho I + \theta T + \ol{\theta} T^* + F_0 = 0$ implies that 
\[
0 = |\theta|^2 i (A - A^*) + F_0. \]
In particular, $A - A^*$ has rank at most $2 m$.   Again, we write $A = R + i B$ where $R = (A+A^*)/2$ and $B = (A-A^*)/2i$.   Then $B$ has rank at most $\mathrm{rank}\, F_0 \le 2m$ and 
\[
T = \ol{\theta} i(R+i B) + \ol{\theta} i \kappa I \]
expresses $T$ in the desired form.
\end{pf}


\begin{lem}\label{lem5.06E}
Let $\hilb$ be a complex Hilbert space and suppose that $T \in \bofh$ satisfies 
\[
\mathrm{rank}\, (\alpha I + \beta T + \gamma T^* + \delta T^* T) \le 3 \]
for some $\alpha, \beta, \gamma, \delta \in \bbC$, where $\delta \ne 0$.   
Then there exist a unitary operator $V$, a finite-rank operator $L$ of rank at most $12$, and $\mu, \lambda \in \bbC$ such that 
\[
T = \lambda (V+ L) + \mu I. \]
\end{lem}

\begin{pf}
It is clear that there is no loss of generality in assuming that $\delta  = \frac{1}{2}$.   Choose $F \in \bofh$ with $\mathrm{rank}\, F \le 3$ such that 
\[
\alpha I + \beta T + \gamma T^* + \frac{1}{2} T^* T + F = 0.\]
This trivially implies that $\ol{\alpha} I + \ol{\gamma} T + \ol{\beta} T^* + \frac{1}{2} T^* T + F^* = 0$.
As before, we set $\rho = \alpha+ \ol{\alpha}$, $\theta = \beta + \ol{\gamma}$ and $F_0 = F + F^*$.    
Then
\[
\rho I + \theta T + \ol{\theta} T  + T^* T + F_0 = 0. \]

A routine calculation shows that 
\[
(\rho - |\theta|^2) I + (T + \ol{\theta}I)^*(T+ \ol{\theta} I) + F_0 = 0. \]
Of course, $(T + \ol{\theta}I)^*(T+ \ol{\theta} I) \ge 0$, and so - by considering this equation modulo the compact operators, we conclude that $|\theta|^2 - \rho \ge 0$.

We again consider two cases.

\bigskip
\noindent{\textsc{Case 1.}} \ \ \ $| \theta|^2 = \rho$.

Then 
\[
(T + \ol{\theta}I)^*(T+ \ol{\theta} I)  = -F_0. \]
But then $|T + \ol{\theta} I|$ has finite rank, so that $G= T + \ol{\theta} I$ has finite rank.   That is, 
\[
T = -\ol{\theta} (I - \ol{\theta}^{-1} G) + 0 I \]
expresses $T$ in the desired form.

\bigskip
\noindent{\textsc{Case 2.}} \ \ \ $| \theta|^2 > \rho$.

Set $V_0 = (|\theta|^2 - \rho)^{-1/2} (T + \ol{\theta} I)$ and $F_2= (|\theta|^2 - \rho)^{-1} F_0$.   Then 
\begin{align*}
V_0^* V_0
	&= (|\theta|^2 - \rho)^{-1} (T+\ol{\theta} I)^* (T + \ol{\theta}I) \\
	&= (|\theta|^2 - \rho)^{-1} ( (|\theta|^2-\rho)I - F_0) \\
	&= I - F_2.
\end{align*}	

That is, $\pi(V_0)$ is an isometry in the Calkin algebra.   In particular, $V_0$ is semi-Fredholm, and therefore $(T+\ol{\theta} I)$ is semi-Fredholm.   But $T + \ol{\theta} I$ has {\offdr}, so $T + \ol{\theta} I$ is biquasitriangular, by Proposition~\ref{prop5.02}.

Hence $V_0$ is Fredholm with index 0.   Using the polar decomposition and the fact that $V_0$ has index 0, we may find a unitary operator $U$ such that  $V_0 = U |V_0| = U (I-F_2)^{1/2}$.   Thus $U - V_0$ is of finite rank, and 
\[
V_0 = U - (U - V_0) = (|\theta|^2 - \rho)^{-1/2} (T + \ol{\theta} I). \]
In other words, 
\[
T = (|\theta|^2 - \rho)^{1/2} (U + (V_0-U)) - \ol{\theta} I \]
again expresses $T$ in the desired form.
\end{pf}


\begin{prop} \label{prop5.07E}
Let $\hilb$ be an infinite-dimensional complex Hilbert space and $F \in \bofh$ be a finite-rank operator.   
\begin{enumerate}
	\item[(a)]
	Suppose that $V \in \bofh$ is unitary. If $W = V + F$ has {\offdr}, then $W$ is normal and $\sigma(W)$ is circlinear.
	\item[(b)]
	Suppose that $R \in \bofh$ is hermitian.  If $L = R+F$ has {\offdr}, then $R$ is normal and $\sigma(L)$ is circlinear.
\end{enumerate}
\end{prop}

\begin{pf}
\begin{enumerate}
	\item[(a)]
	It is obvious that $I - W^* W$ is of finite rank.  By Corollary~6.17 of~\cite{RadR2003}, $W$ must have a non-trivial invariant subspace, which  must -- by virtue of {\offdr} --  in fact be an orthogonally reducing subspace for $W$.  Thus we may write $W \simeq W_1 \oplus W_2$, and it is clear that each $W_k$ must satisfy {\offdr} (by Proposition~\ref{prop2.02}) and be of the form $V_k + F_k$ for some unitary operator $V_k$ and some finite-rank operator $F_k$, $k = 1,2$.  At least one of these summands acts on an infinite-dimensional space, and thus we may again apply Theorem~6.17 of~\cite{RadR2003} to find non-trivial invariant -- hence reducing -- subspaces for that summand.
	
	Repeating this process, we see that for any $n \ge 1$, we can find $n$ summands $X_{n, 1},  X_{n,2}, \ldots,$ $X_{n, n}$ of $W$ such that 
	\[
	W = X_{n,1} \oplus X_{n,2} \oplus \cdots \oplus X_{n,n}. \]
	Furthermore, a moments' thought will convince the reader that at most $\mathrm{rank}\, F$ of these summands can fail to be unitary themselves, and hence when $n > \mathrm{rank}\, F$, at least one of the $X_{n, k}$'s is a unitary operator.
	
	Let \[
	\ \ \ \ \ \ \ \ \cJ = \{ (U, \cM): U \mbox{ is a unitary direct summand of } W \mbox{ acting on the subspace } \cM \mbox{ of } \hilb\}. \]
	The above paragraph shows that $\cJ$ is non-empty.  Partially order $\cJ$ by setting $(U_1, \cM_1) \le (U_2, \cM_2)$ if $\cM_1 \le \cM_2$.  (Note that this automatically implies that $U_1$ is a direct summand of $U_2$.)  If $\cC = \{ (U_\nu, \cM_\nu) : \nu \in \Gamma \}$ is a chain in $\cC$, then by setting $\cM = \ol {\cup_{\nu \in \Gamma} \cM_\nu}$, we see that $\cM$ is a reducing subspace for $W$ (as each $\cM_\nu$ is), and $U = W|_{\cM}$ is unitary (since it is clearly unitary on the dense submanifold $\cup_{\nu \in \Gamma} \cM_\nu$ of $\cM$).  It follows from Zorn's Lemma that $\cJ$ admits a maximal element $(U_0, \cM_0)$.   If $\cM_0^\perp$ is infinite-dimensional, then the argument of the first two paragraphs can be used to show that $W|_{\cM_0^\perp}$ admits a unitary direct summand, contradicting the maximality of $(U_0, \cM_0)$.  	Thus $m= \dim\, \cM_0^\perp < \infty$. 
	
	Write $W = U_0 \oplus Y$, where $Y$ acts on $\cM_0^\perp$, and note that $Y$ has {\offdr}.  We may view $Y$ as an element of $\bbM_m(\bbC)$, so that $Y$ can be upper triangularized with respect to some orthonormal basis.  The fact that $Y$ has {\offdr} implies that it is reductive, and is therefore normal.   This forces $W$ to be normal as well.   There remains to show that $\sigma(W)$ is circlinear.  Note that if $\sigma(W)$ is finite, then all elements of $\sigma(W)$ are eigenvalues, and so $\sigma(W)$ is circlinear by the comments of Section~\ref{sec5.03}.
	
	Hence we may assume that $\sigma(W)$ is infinite, which is equivalent to assuming that $\sigma(U_0)$ is infinite.   In this case, we shall prove that $W$ is unitary.  We argue by contradiction.  Suppose otherwise, and let $\tau \in \sigma(Y)$ with $|\tau| \ne 1$.  Let $\cN \subseteq \ker\, (W - \tau I) \subseteq \cM_0^\perp$ be a one-dimensional subspace.  We see that the operator $Z = U_0 \oplus \tau$, being a direct summand of $W$, also satisfies {\offdr}.   With respect to the decomposition $\cM_0 \oplus \cN$, we may write 
	\[
	Z = \begin{bmatrix} U_0 & 0 \\ 0 & \tau\end{bmatrix}. \]
	
	Let $x \in \cM_0$ be a unit vector such that $\{ x, U_0 x, U_0^2 x\}$ is linearly independent.   Such a vector must exist, otherwise $U_0$ is \emph{boundedly locally linearly dependent}, which -- by Kaplansky's Theorem~\cite{Kap1971}, Lemma~14 -- implies that $U_0$ is algebraic, and therefore has finite spectrum, a contradiction of our current assumption.
	
	Thus $\{ U_0^* x, x, U_0 x \}$ is again linearly independent, as $U_0$ is unitary.   We shall now find vectors $y$ and $z$ in $\cM_0 \oplus \cN$ such that $\cE_1 = \{ y, z, Zy, Zz \}$ is linearly independent, but $\cE_2 = \{ y, z, Z^*y, Z^* z \}$ is not.   This will yield the desired contradiction, by implying that $(I-P)Z P$ and $(I-P)Z^* P$ have ranks two and one respectively.  
	
	Let $y = \begin{bmatrix} x \\ 1 \end{bmatrix}$ and $z = \begin{bmatrix} U_0 x \\\xi \end{bmatrix}$, with $\xi \in \bbC$ to be determined shortly.  (Here, we have identified $\cN$ with $\bbC$.)  Now 
	\[
	\cE_1 = \left \{ \begin{bmatrix} x \\ 1 \end{bmatrix}, \begin{bmatrix} U_0 x \\\xi \end{bmatrix}, \begin{bmatrix} U_0 x \\ \tau \end{bmatrix}, \begin{bmatrix}  U_0^2 x \\ \tau \xi \end{bmatrix} \right\}, \]
	and
	\[
	\cE_2 = \left \{ \begin{bmatrix} x \\ 1 \end{bmatrix}, \begin{bmatrix}  U_0 x \\\xi \end{bmatrix}, \begin{bmatrix}  U_0^* x \\ \ol{\tau} \end{bmatrix}, \begin{bmatrix} x \\ \ol{\tau} \xi \end{bmatrix} \right\}. \]

Let $\xi = \tau$.   Then $\cE_1$ is linearly dependent, but $\cE_2$ is not, because $\ol{\tau} \xi = |\tau|^2 \ne 1$.
	\item[(b)]
	The proof of this result is similar.   We may use Corollary~6.15 of~\cite{RadR2003} to assert that if $L - L^*$ has finite rank, then $L$ has a non-trivial invariant subspace, which is again orthogonally reducing by our hypothesis that $L$ satisfies {\offdr}.   One then looks for a maximal hermitian direct summand, and separately argues the cases where that summand has finite or infinite spectrum.  The details are left to the reader.
\end{enumerate}
\end{pf}


\begin{thm} \label{thm5.08E}
Let $\hilb$ be an infinite-dimensional, complex Hilbert space, and let $T \in \bofh$.  If $T$ satisfies {\offdr}, then there exist $\lambda, \mu \in \bbC$ and $A \in \bofh$ with $A$ either selfadjoint or an orthogonally reductive unitary operator such that  $T = \lambda A + \mu I$.

In particular, if $T$ satisfies {\offdr}, then $T$ is normal with circlinear spectrum.
\end{thm}

\begin{pf}
By combining Lemma~\ref{lem5.05E} and Lemma~\ref{lem5.06E}, we can assume without loss of generality that $T = X + F$, where $F$ is of finite rank and $X$ is either selfadjoint or unitary.  
	
	Either way, by Proposition~\ref{prop5.07E}, we see that $T$ is normal with circlinear spectrum.  From this it is easy to verify that $T$ is of the form $\lambda A + \mu I$ for some $\lambda, \mu \in \bbC$ with $A$ either selfadjoint or unitary.   The fact that $T$ is orthogonally reductive implies that $A$ is as well.  (This last argument is superfluous when considering the case where $A$ is selfadjoint.)
\end{pf}


\begin{eg} \label{eg5.09E}
We mention in passing that {\offdr} is a strictly stronger condition than that of being orthogonally reductive.   Indeed, suppose that $N \in \bofh$ is a normal operator with $\sigma(N) = \{ 1, 2, 3, 4+i\}$.  Thus the eigenvalues of $N$ are neither co-linear nor co-circular, and so $N$ does not have {\offdr}, by Corollary~\ref{cor3.15}.   However, then $N$ is orthogonally reductive, as $N^*$ is a polynomial function of $N$, combined with Sarason's result~\cite{Sar1966}.   
\end{eg}


\subsection{} \label{sec5.10E}
It would be interesting to know whether or not the converse of Theorem~\ref{thm5.08E} holds.  

Indeed, suppose that $N \in \bofh$ is normal and has co-linear spectrum.   Arguing as before, we have that there exist scalars $\lambda, \mu \in \bbC$ and a hermitian operator $L$ such that $N = \lambda I + \mu L$.  It is routine to verify that $N$ has {\offdr}.

For normal operators with co-circular spectrum, the problem is a bit more complicated.


\begin{prop} \label{prop5.11E}
Let $U \in \bofh$ be unitary and suppose that $\sigma (U) \ne \bbT$.   Then $U$ has {\offdr}.  
\end{prop}

\begin{pf}  
Let $0 \ne P \ne I$ be a projection in $\bofh$, and relative to $\hilb = P\hilb \oplus P^\perp \hilb$, write 
\[
U = \begin{bmatrix} A & B \\ C & D \end{bmatrix}. \]

Our goal is to show that $\mathrm{rank}\, B = \mathrm{rank}\, C$.   As always, we have 
\begin{align*}
B B^* &= I - A A^* \\
C^* C &= I - A^* A.
\end{align*}

If $B$ and $C$ are both of infinite rank, then there is nothing to prove.  Thus we may suppose that either $B$ or $C$ is of finite rank.    Now, since $U$ has {\offdr} if and only if $U^*$ has {\offdr}, we may suppose -- by taking adjoints if necessary -- that $\mathrm{rank}\, C$ is of finite rank and that $\mathrm{rank}\, C  \le \mathrm{rank}\, B$.  

\bigskip

\noindent{\textbf{Case 1.}}
$B$ is compact.

\smallskip

Then $A$ and $D$ are essentially unitary.   Since $\sigma_e(A) \subseteq \sigma_e(U) \ne \bbT$, it follows that $\mathrm{ind}\, A = 0$.   (That is, in order for $A$ to have non-zero index, $0$ must lie in a bounded component of $\bbC \setminus \sigma_e(A)$, of which there are none.)

Write $A = V |A|$, and note that as $\mathrm{ind}\, A = 0$, we may assume without loss of generality that $V$ is unitary.   Thus $A A^* = V |A| |A| V^* = V (A^* A) V^*$.   That is, $A A^*$ and $A^* A$ are unitarily equivalent.   

But then $B B^*$ and $C^* C$ are unitarily equivalent, whence $\mathrm{rank}\, B = \mathrm{rank}\, B B^* = \mathrm{rank}\, C^* C = \mathrm{rank}\, C$.

\noindent{\textbf{Case 2.}}
$B$ is not compact.

\smallskip
We shall show that under the hypothesis that $\sigma(U) \ne \bbT$, this cannot happen.   Indeed, the equation $C^* C = I - A^*A$ with $\mathrm{rank}\, C < \infty$ implies that 
\[
1 = \pi(I) = \pi(A)^* \pi(A). \]
Thus $\pi(A)$ is a partial isometry in the Calkin algebra $\bofh/\kofh$, which implies that $\pi(A) \pi(A)^*$ is a projection.
The fact that $B$ is not compact, combined with the fact that $B B^* = I - A A^*$ shows that 
\[
1 = \pi(I) \ne \pi(A) \pi(A)^*. \]
Thus $\pi(A)$ is not unitary.
Choose a projection $R \in \bofh$ such that $\pi(R) = \pi(A) \pi(A)^*$.   By Lemma~V.6.4 of~\cite{Dav1996}, there exists a partial isometry $W \in \bofh$ such that $W = R W$ and $\pi(W)=\pi(A)$.  Moreover, by that same result, the integer 
\[
\xi =\mathrm{rank}\, (I - W^* W) - \mathrm{rank}\, (R - W W^*)\]
is defined independent of the choice of $W$.

In our case, $\mathrm{rank}\, (I - W^* W) < \infty$ while $\mathrm{rank}\, (I-R) = \infty$ and $\mathrm{rank}\, (R- W W^*) < \infty$.  Hence $\mathrm{rank}\, (I - W W^*) = \infty$.

Thus we have that $W$ is a partial isometry with initial space $W^* W \hilb$, and final space $W \hilb$, and 
\begin{enumerate}
	\item[(i)]
	$\dim\, (W^* W \hilb)^\perp < \infty$; and 
	\item[(ii)]
	$\dim\, (W \hilb)^\perp = \infty$.   
\end{enumerate}
It is routine to produce a partial isometry $W_0$ with initial space $(W^* W \hilb)^\perp < \infty$ and final space contained in $(W \hilb)^\perp$, and to verify that $V = W + W_0$ is an isometry on $\hilb$.

By the Wold Decomposition, $V$ is unitarily equivalent to $S^{(\kappa)} \oplus Y$, where $S$ denotes the unilateral forward shift, $Y$ is a unitary operator, and $\kappa \in \bbN \cup \{ 0, \infty\}$.  

If $\kappa = 0$, then $V$ is unitary.  But then $\pi(V)= \pi(W)= \pi(A)$ is also unitary, a contradiction.   Thus $\kappa \ne 0$.   But then $\sigma_e(A) = \sigma_e(V) \supseteq \sigma_e(S) = \bbT$.

On the other hand, it is not too hard to show that $\partial(\sigma_e(A)) \subseteq \sigma_{\ell r e}(A) \subseteq \sigma_e(U)$. (For example, by the Corollary to Theorem~4.3 of~\cite{FSW1972}, there exists a compact operator $K_1 \in \cB(P \hilb)$ such that $A + K_1 = \begin{bmatrix} \lambda I & 0 \\ 0 & A_4 \end{bmatrix}$ with respect to the decomposition $P \hilb = \cM \oplus (P\hilb \ominus \cM)$, for an appropriate subspace $\cM \subseteq P \hilb$ satisfying $\dim\, \cM = \dim\, (P \hilb \ominus \cM) = \infty$.   Letting $K = K_1 \oplus 0$ yields that $U + K = \begin{bmatrix} \lambda I & 0 & B_1 \\ 0 & A_4 & B_2 \\ 0 & 0 & D \end{bmatrix}$.   Thus $\lambda \in \sigma_e(U+K) = \sigma_e(K)$.) But  $\partial(\sigma_e(A)) = \bbT$, which contradicts our hypothesis that $\sigma(U) \ne \bbT$.

This shows that the case where $C$ is of finite rank and $B$ is not compact cannot happen, and completes the proof.
\end{pf}


Having seen that the bilateral shift $W$ is a unitary operator with $\sigma(W) = \bbT$ which is not reductive, we now show that there exists a unitary operator whose spectrum is the unit circle $\bbT$, but which nonetheless has  {\offdr}.

Before embarking upon the proof of this, we first require a result due to Wu and Takahashi~\cite{WT1999}.  Recall that if $X \in \bofh$ is an operator, then we define the \textbf{defect indices} of $X$ to be
\begin{align*}
	d_X &= \dim (\ol{\mathrm{ran}}(I-X^* X)^{1/2}) \mbox{ and } \\
	d_{X^*} &= \dim (\ol{\mathrm{ran}} (I-X X^*)^{1/2}).
\end{align*}


\begin{prop} \label{5.12E} \emph{(}\textbf{Wu-Takahashi}; Theorem~3.5 of~\cite{WT1999}\emph{)}
Let $X \in \bofh$ be a contraction and suppose that $d_X \ne d_{X^*}$.  Then $X$ does not admit a singular unitary dilation.
\end{prop}


\begin{prop} \label{prop5.13E}
Let $ (d_n)_{n=1}^\infty$ be a sequence in $\bbT$ and let $V = \mathrm{diag}\, (d_n)_{n=1}^\infty$ be a corresponding diagonal unitary operator in $\bofh$.  Then $V$ has {\offdr}.
\end{prop}

\begin{pf}
Suppose that we can find a projection $P \in \bofh$ such that with respect to the decomposition $\hilb = P \hilb \oplus P^\perp \hilb$ we may write
\[
V = \begin{bmatrix} A & B \\ C &D \end{bmatrix} \]
with $\mathrm{rank}\, C < \mathrm{rank}\, B$.   

(In particular, we must have $\mathrm{rank}\, C < \infty$).   Then 
\begin{align*}
	B B^* &= I - A A^* \mbox{ and }\\
	C^* C &= I - A^* A
\end{align*}	
have different ranks.  Since $C$ is of finite rank, $A$ is essentially isometric and thus is a semi-Fredholm operator - in particular, both $A$ and $A^*$ have closed range.   Also, 
\begin{align*}
	\mathrm{rank}\, C = \mathrm{rank}\, C^* C &= \mathrm{rank}\, (I - A A^*) = d_{A^*} < \infty, \mbox{ while } \\
	\mathrm{rank}\, B = \mathrm{rank}\, B B^* &= \mathrm{rank}\, (I- A^* A) = d_{A}.
\end{align*}	

In other words, $A$ is a contraction and the defect indices of $A$ are unequal.   By 
Proposition~\ref{5.12E} above, $A$ does not admit unitary dilation.   But $U$ is diagonal, and is therefore a singular  unitary dilation of $A$, which is obviously a contradiction.
\end{pf}


If, in Proposition~\ref{prop5.13E} we choose $\{ d_n\}_n$ to be dense in $\bbT$, we immediately obtain the following consequence:

\begin{cor} \label{cor5.14E}
There exists a unitary operator $V$ with $\sigma(V)= \bbT$ which has {\offdr}.
\end{cor}


\subsection{} \label{sec5.15E}
Wermer~\cite{Wer1952} has shown that a unitary operator $U$ fails to be reductive if and only if Lebesgue measure is absolutely continuous with respect to the spectral measure $\mu$ for $U$.  Since any operator with {\offdr} is necessarily reductive, this provides a measure-theoretic obstruction to {\offdr} for unitary operators.

\smallskip

Another consequence of the above analysis is that it proves that the set $\mathfrak{G}_{\mathfrak{rank}}$ of operators with {\offdr} is not closed.  Indeed, it follows easily from~\cite{Dav1986} that the bilateral shift $W$ is a limit of unitary operators $V_n$ such that $\sigma(V_n) \ne \bbT$.  (The $V_n$'s can in fact be chosen to be  unitary operators with spectrum  $\Gamma_n = \{ e^{2 \pi i \theta}: 0 \le \theta \le 1 - \frac{1}{n} \}$.)   As we saw in Proposition~\ref{prop5.11E}, each $V_n$ has {\offdr}, but $W = \lim_n V_n$ does not. 

Alternatively, the Weyl-von Neumann-Berg Theorem (see, e.g.,~\cite{Dav1996}, Theorem~II.4.4) shows that there exists a sequence $(W_n)_{n=1}^\infty$ of diagonal unitary operators such that $\sigma(W_n) = \bbT$ for all $n \ge 1$, such that $W = \lim_n W_n$.

\bigskip

We now investigate a consequence of {\offdr} which relates to cyclic subspaces for operators.


\begin{prop}\label{prop5.16E}
Suppose that $T$ is reductive.  Then $T$ and $T^*$ have the same cyclic subspaces.    In particular, if $T$ has {\offdr}, then $T$ and $T^*$ have the same cyclic subspaces.
\end{prop}

\begin{pf}
Suppose that $0 \ne \cM \subseteq \hilb$ is a cyclic subspace for $T$, and let $0 \ne x \in \cM$ be a cyclic vector for $T$ in $\cM$, so that $\cM = \ol{\mathrm{span}} \{ x, T x, T^2 x, \ldots \}$.   If $P$ is the orthogonal projection of $\hilb$ onto $\cM$, then $P^\perp T P = 0$, so by reductivity, $P T P^\perp = 0$, which implies that $\cM$ is invariant for $T^*$.

Now let $\cN = \ol{\mathrm{span}}\, \{ x, T^*x, (T^*)^2 x, \ldots \}$ be the cyclic subspace for $T^*$ generated by $x$.  Since $x \in \cM$ and $\cM$ is invariant for $T^*$, we see that $\cN \subseteq \cM$.  Also, as $T^*$ is also reductive, the argument of the first paragraph shows that $x \in \cN$ is invariant for $T$.   But then $\cN \supseteq \cM$, whence $\cN = \cM$, completing the proof.  
\end{pf}


There exists a variant of this result which is somewhat interesting.

\begin{prop}  \label{prop5.17E}
Let $T \in \bofh$ and suppose that for each orthogonal projection $P \in \bofh$, the off-diagonal corner $P^\perp T P$ has rank one if and only if $P T P^\perp$ has rank one.   A subspace $\cM$ of $\hilb$ of dimension at least 3 is cyclic for $T$ if and only if it is cyclic for $T^*$.  
\end{prop}

\begin{pf}
Given $T$ as in the statement of the Proposition, it is clear that $T^*$ also has this property.   The argument used to prove Proposition~\ref{prop5.16E} shows that it suffices to show that the cyclic subspace $\cM = \ol{\mathrm{span}} \{ x, Tx, T^2 x, \ldots \}$ for $T$ generated by a non-zero vector $x$ is invariant for $T^*$. 

We consider first the case where $\cM$ is infinite-dimensional, as it is the easier of the two.

For each $n\ge 1$, let $P_n$ denote the orthogonal projection of $\hilb$ onto $\ol{\mathrm{span}} \{ x, Tx, T^2 x, \ldots, T^n x \}$, and let $P$ denote the orthogonal projection of $\hilb$ onto $\cM$.   It is clear that $(P_n)_{n=1}^\infty$ is an increasing sequence which converges in the strong operator topology to $P$.      An easy calculation then shows that the sequence $(P_n^\perp T P_n)_{n=1}^\infty$  converges in the strong operator topology to $P^\perp T P$.    

As $x$ is a cyclic vector for $T$ in $\cM$, we have that $\mathrm{rank}\, (P_n^\perp T P_n) = 1$ for all $n \ge 1$, and our hypothesis then asserts that $\mathrm{rank}\, (P_n^\perp T^* P_n) = 1$ for all $n \ge 1$.  But rank is lower-semicontinuous with respect to the strong operator topology, and thus $\mathrm{rank}\, P^\perp T^* P \le 1$.   If $\mathrm{rank}\, P^\perp T^* P = 1$, then the hypothesis on $T$ implies that $\mathrm{rank}\, P^\perp T P = 1$, contradicting the fact that $\cM$ is invariant for $T$.   Hence $P^\perp T^* P = 0$, proving that $\cM$ is invariant for $T^*$.  

\smallskip

Next we suppose that $\cM$ is finite-dimensional with $\dim\, \cM = N \ge 3$, and we find a cyclic vector $x$ for $T$ so that  $\cM =\mathrm{span} \{ x, Tx, T^2 x, \ldots, T^{N-1} x\}$.   Let $\{ e_1, e_2, \ldots, e_N\}$ be the orthonormal basis obtained from $\{ x, Tx, T^2 x, \ldots, T^{N-1} x\}$ by applying the Gram-Schmidt process, so that $\mathrm{span} \{ e_1, e_2, \ldots, e_k\}= \mathrm{span} \{x, Tx, \ldots, T^{k-1} x \}$ for $1 \le k \le N$.  Let $Q_k$ denote the orthogonal projection of $\hilb$ onto $\bbC e_k$, $1 \le k$, and define $P_k = Q_1 + Q_2 + \cdots + Q_k$, $1 \le k \le N$.  Finally, extend $\{ e_k\}_{k=1}^N$ to an orthonormal basis $\{ e_k \}_{k=1}^\infty$ for $\hilb$.

Note that the fact that $x$ is cyclic for $\cM$, combined with our hypothesis,  implies that $\mathrm{rank}\, P_k^\perp T P_k = 1 = \mathrm{rank}\, P_k T P_k^\perp$, $1 \le k \le N-1$.  Moreover, $\cM$ is invariant for $T$, whence $P_N^\perp T P_N = P^\perp T P = 0$.   By hypothesis, $\mathrm{rank}\, P T P^\perp \ne 1$.   

But 
\begin{align*}
P T P^\perp 
	&= P_N T P_N^\perp \\
	&= P_{N-1} T P_{N-1}^\perp P_N^\perp + Q_N T P_N^\perp,
\end{align*}
so that $\mathrm{rank}\, P T P^\perp \le \mathrm{rank}\, P_{N-1} T P_{N-1}^\perp + \mathrm{rank}\, Q_N T P_N^\perp \le 1 + 1 = 2.$

Thus $\mathrm{rank}\, P T P^\perp \in \{ 0, 2\}$, and our goal is to show that $\mathrm{rank}\, P T P^\perp \ne 2$.

\smallskip

Suppose, to the contrary, that $\mathrm{rank}\, P T P^\perp = 2$.  It follows that $\mathrm{rank}\, P_{N-1} T P_{N-1}^\perp = 1 = \mathrm{rank}\, Q_N T P_N^\perp$.   Thus there exists $1 \le k \le N-1$ such that $Q_k T P_N^\perp \ne 0$, and $\mathrm{rank}\, (Q_k + Q_N) T P^\perp =\mathrm{rank}\, (Q_k + Q_N) T P_N^\perp = 2$.   

\smallskip

\noindent{\textsc{Case 1.}} $k = N-1$.
Let us reorder the basis for $P \hilb$ as $\{ e_{N-1}, e_N, e_1, e_2, \ldots, e_{N-2} \}$.   The matrix for $T$ relative to $ P \hilb \oplus P^\perp \hilb$ is:
\[
[T] = \begin{bmatrix} t_{N-1, N-1} & t_{N-1, N} & t_{N-1, 1} & t_{N-1, 2} & \ldots & t_{N-1, N-3} & t_{N-1, N-2} & Q_{N-1} T P^\perp \\
t_{N, N-1} & t_{N, N} & t_{N, 1} & t_{N, 2} & \ldots & t_{N, N-3} & t_{N, N-2} & Q_{N} T P^\perp \\
\vdots &   &  &  &\ldots  &   & \vdots \\
t_{N-3, N-1} & t_{N-3, N} & t_{N-3, 1} & t_{N-3, 2} & \ldots & t_{N-3, N-3} & t_{N-3, N-2} & Q_{N-3} T P^\perp \\
t_{N-2, N-1} & t_{N-2, N} & t_{N-2, 1} & t_{N-2, 2} & \ldots & t_{N-2, N-3} & t_{N-2, N-2} & Q_{N-2} T P^\perp \\
0 & 0 & 0 & 0 & \ldots & 0 & 0 & P^\perp T P^\perp
\end{bmatrix}\]

Let $R = P - Q_{N-2}$.  Since $t_{N-2, N-3} \ne 0$ (as $x$ is a cyclic vector for $\cM$), it follows that $\mathrm{rank}\, R^\perp T R = 1$.

Thus \[
\mathrm{rank}\, \begin{bmatrix} t_{N-1, N-2} & Q_{N-1} T P^\perp \\ t_{N, N-2} & Q_N T P^\perp \\ \vdots & \vdots \\ t_{N-3, N-2} & Q_{N-3} T P^\perp \end{bmatrix} = \mathrm{rank}\, R T R^\perp = 1, \]
 and so $\mathrm{rank}\, (Q_{N-1} + Q_N) T P^\perp = \mathrm{rank} \begin{bmatrix} Q_{N-1} T P^\perp \\ Q_N T P^\perp \end{bmatrix} \le 1$, a contradiction.  Thus in this case, $P T P^\perp = 0$, so $\cM$ is invariant for $T^*$.
\bigskip


\noindent{\textsc{Case 2.}} $1 \le k < N-1$.

This time we reorder the basis for $P \hilb$ as $\{ e_k, e_{k+2}, \ldots,  e_{N}, e_1, \ldots, e_{k-1}, e_{k+1} \}$. The matrix for $T$ relative to $P \hilb \oplus P^\perp \hilb$ is then:
\[
[T] = \begin{bmatrix} t_{k,k}  &  t_{k, k+2} & \ldots & t_{k, N}  & t_{k, 1} & \ldots & t_{k, k-1} &t_{k, k+1}  & Q_{k} T P_N^\perp \\
t_{k+2,k}  & t_{k+2, k+2} & \ldots & t_{k+2, N}  & t_{k+2, 1} & \ldots & t_{k+2, k-1} &t_{k+2, k+1}  & Q_{k+2} T P_N^\perp \\
\vdots &   &  &  &\ldots  &   & \vdots \\
t_{k-1,k}  & t_{k-1, k+2} & \ldots & t_{k-1, N}  & t_{k-1, 1} & \ldots & t_{k-1, k-1} &t_{k-1, k+1}  & Q_{k-1} T P_N^\perp \\
t_{k+1,k}  & t_{k+1, k+2} & \ldots & t_{k+1, N}  & t_{k+1, 1} & \ldots & t_{k+1, k-1} &t_{k+1, k+1}  & Q_{k+1} T P_N^\perp \\
0 & 0 & 0 & 0 & 0 & \ldots & 0 & 0 & P_N^\perp T P_N^\perp
\end{bmatrix}\]

Let $R = P - Q_{k+1}$.  Since $t_{k+1, k} \ne 0$ (as $x$ is a cyclic vector for $\cM$), it follows that $\mathrm{rank}\, R^\perp T R = 1$.  By hypothesis, $\mathrm{rank}\, R T R^\perp = 1$.

Thus 
\[
\mathrm{rank}\, (Q_k +Q_N) T P^\perp 
	=\mathrm{rank}\, (Q_k + Q_N) [ R T  R^\perp] P^\perp \le \mathrm{rank}\, R T R^\perp = 1, \]
a contradiction.  	  Thus in this case as well, $P T P^\perp = 0$, so $\cM$ is invariant for $T^*$.

\bigskip

The remainder of the proof is identical to the second paragraph of the proof of Proposition~\ref{prop5.16E}.
\end{pf}
		


\section{Essentially reductive operators with {\offdr}} \label{sec6}


\subsection{} \label{sec6.01}
In~\cite{Har1975}, Ken Harrison introduced the notion of \textbf{essentially reductive} operators:  we say that $T \in \bofh$ is \textbf{essentially reductive} if for each projection $P$ we have that $P T P^\perp$ compact if and only if $P^\perp T P$ is compact.   (One can view this as $\pi(T)$ having {\offdr} in the Calkin algebra.)

In the paper~\cite{Moo1978} (Theorem~2), Moore shows that every essentially reductive operator $T$ is essentially normal -- i.e. $\pi(T)$ is normal in the Calkin algebra.  Earlier, Harrison (\cite{Har1975}, Theorem~4.5) had characterized all essentially normal operators which are essentially reductive.   Combining these results, one obtains the following.


\begin{thm} \label{thm6.02} \emph{(}\textbf{Moore}; Corollary~1 of~\cite{Moo1978}\emph{)}
Let $T \in \bofh$.   The following are equivalent.
\begin{enumerate}
	\item[(a)]
	$T$ is essentially reductive.
	\item[(b)]
	$T$ is essentially normal and $\sigma_e(T)$ is Lavrentiev.
\end{enumerate}
\end{thm}	


The next result is a simple consequence of Moore's Theorem together with Theorem~\ref{thm5.08E} and Proposition~\ref{prop5.11E}.

\begin{cor} \label{cor6.03}
Let $T \in \bofh$.  The following are equivalent.
\begin{enumerate}
	\item[(a)]
	$T$ is essentially reductive and has {\offdr}.
 	\item[(b)]
	One of the following holds.
		\begin{enumerate}
		\item[(i)]
		There exist $\lambda, \mu \in \bbC$ and a hermitian operator $R$ such that $T = \lambda R  + \mu I$. 
		\item[(ii)]
		There exist $\lambda, \mu \in \bbC$ and a unitary operator $V$ with $\sigma(V) \ne \bbT$ such that $T = \lambda V + \mu I$.
		\end{enumerate}
\end{enumerate}
\end{cor}
		




\providecommand{\bysame}{\leavevmode\hbox to3em{\hrulefill}\thinspace}
\providecommand{\MR}{\relax\ifhmode\unskip\space\fi MR }
\providecommand{\MRhref}[2]{%
  \href{http://www.ams.org/mathscinet-getitem?mr=#1}{#2}
}
\providecommand{\href}[2]{#2}

\end{document}